\input amstex\documentstyle{amsppt}  
\pagewidth{12.5cm}\pageheight{19cm}\magnification\magstep1
\topmatter
\title Total positivity in reductive groups, II \endtitle
\author G. Lusztig\endauthor
\address{Department of Mathematics, M.I.T., Cambridge, MA 02139}\endaddress
\thanks{Supported by NSF grant DMS-1566618.}\endthanks
\endtopmatter   
\document

\define\uI{\un I}

\define\ui{\un i}
\define\uj{\un j}

\define\dw{\dot w}

\define\ds{\dot s}

\define\dW{\dot W}
\define\dI{\dot I}

\define\mpb{\medpagebreak}

\define\frl{\forall}

\define\si{\sim}

\define\sqc{\sqcup}

\define\qua{\quad}

\define\lb{\linebreak}

\define\op{\oplus}
   
\define\part{\partial}
\define\emp{\emptyset}

\define\ra{\rangle}
\define\n{\notin}
\define\iy{\infty}
\define\m{\mapsto}
\define\do{\dots}
\define\la{\langle}

\define\sm{\smallmatrix}
\define\esm{\endsmallmatrix}
\define\sub{\subset}    

\define\T{\times}
\define\ti{\tilde}
\define\nl{\newline}
\redefine\i{^{-1}}
\define\fra{\frac}
\define\un{\underline}

\define\ot{\otimes}

\define\a{\alpha}
\redefine\b{\beta}
\redefine\c{\chi}
\define\g{\gamma}
\redefine\d{\delta}
\define\e{\epsilon}
\define\et{\eta}
\define\io{\iota}
\redefine\o{\omega}
\define\p{\pi}
\define\ph{\phi}
\define\ps{\psi}
\define\r{\rho}
\define\s{\sigma}
\redefine\t{\tau}
\define\th{\theta}
\define\k{\kappa}
\redefine\l{\lambda}
\define\z{\zeta}
\define\x{\xi}

\redefine\G{\Gamma}

\redefine\L{\Lambda}
\define\Ph{\Phi}
\define\Ps{\Psi}

\redefine\aa{\bold a}

\define\hh{\bold h}
\define\ii{\bold i}
\define\jj{\bold j}
\define\kk{\bold k}

\define\BB{\bold B}
\define\CC{\bold C}

\define\NN{\bold N}

\define\QQ{\bold Q}
\define\RR{\bold R}

\define\UU{\bold U}

\define\ZZ{\bold Z}

\define\cb{\Cal B}
\define\cc{\Cal C}

\define\cj{\Cal J}

\define\cn{\Cal N}
\define\co{\Cal O}
\define\cp{\Cal P}

\define\cv{\Cal V}

\define\cz{\Cal Z}
\define\cx{\Cal X}
\define\cy{\Cal Y}

\define\ff{\frak f}

\define\fC{\frak C}

\define\fG{\frak G}

\define\fP{\frak P}

\define\fU{\frak U}

\define\ta{\ti a}

\define\tf{\ti f}

\define\tx{\ti x}

\define\tB{\ti B}

\define\tS{\ti S}

\define\sha{\sharp}

\define\cir{\circ}

\head Introduction \endhead
\subhead 0.1\endsubhead
The theory of totally positive real matrices (see \cite{P11}) has been developed in the
1930's by I.Schoenberg and independently by F.Gantmacher and M.Krein after earlier contributions by
M.Fekete and G.Polya (1912) and with later contributions by A.Whitney and C.Loewner (1950's). 
In \cite{L94} I extended a part of this theory by replacing the group $SL_n(\RR)$ by an arbitrary
split semisimple real Lie group $G$. This paper is a continuation of \cite{L94}.

\subhead 0.2\endsubhead
One of the main tools in \cite{L94} was the use of the canonical basis $\BB$ of the $+$ part $\UU^+$ of a quantized enveloping algebra $\UU$ 
(over $\QQ(v)$) of type $A,D,E$, introduced in \cite{L90}, and in particular the use of its positivity properties.
Now $\UU^+$ admits a family of bases (the PBW-bases) indexed by the various reduced expressions for the longest element $w_I$ in 
the Weyl group. In \cite{L90} it was shown that the specialization of any of these PBW-bases at $v=\iy$
is independent of the PBW-basis and is in fact the same as the specialization of $\BB$ at $v=\iy$. Since each 
PBW-basis is naturally parametrized by $\NN^\nu$ (where $\nu$ is the number of positive roots), it follows that 
$\BB$ has a family of parametrizations by $\NN^\nu$ indexed by the various reduced expressions for $w_I$. 
Moreover, in \cite{L90} it is shown that any two of these parametrizations are related by a piecewise linear 
automorphism of $\NN^\nu$ built from adding or substracting two integers or by taking the minimum of two integers. 
These piecewise linear automorphisms extend in an obvious way to piecewise linear automorphisms of $\ZZ^\nu$ which
can be viewed as maps relating the parametrizations of a new object $\ti\BB$ (containing $\BB$ as a subset)
by $\ZZ^\nu$, these parametrizations being again indexed by the various reduced expressions for $w_I$. 
The set $\ti\BB$ with this family of parametrizations by $\ZZ^\nu$ is the first example of a positive structure 
(see 1.3), which is one of the themes of this paper. 

A second example of positive structure appeared in \cite{L94}. In \cite{L94} I define the positive part 
$U^+_{>0}$ and the non-negative part $U^+_{\ge0}$ of the strictly upper triangular part $U^+$ of the real Lie 
group $G$ (which in this introduction is assumed for simplicity to be of the same type as $\UU$ above). In 
\cite{L94} it is shown that $U^+_{>0}$ has several parametrizations by $\RR_{>0}^\nu$; moreover any two of these 
parametrizations are related by an automorphism of $\RR_{>0}^\nu$ which is built from multiplying or dividing 
two numbers in $\RR_{>0}$, or adding two numbers in $\RR_{>0}$. Moreover, in \cite{L94} it is shown that this 
last automorphism is related to the analogous automorphism for $\ti\BB$ by a process (``passage to zones'')
 which connects geometrical objects over $\RR(t)$ ($t$ is an indeterminate) with piecewise linear objects 
involving only integers. I believe that this was the first time that such a connection was used in relation to 
Lie theory, see also 1.4. In later works, this process, introduced in \cite{L94}, appears under the name
of ``tropicalization'' (see for example \cite{FG}) or under the name of ``ultra-discretization'' (see for 
example \cite{NY}).

We now describe the contents of the various sections.

In \S1 we give a definition of a positive $K$-structure and that of a non-negative $K$-structure. Here $K$ can be
for example $\RR_{>0}$ or $\ZZ$, see 1.1. 

In \S2 we recall following \cite{L94} the definition of the non-negative
submonoid of the real Lie groups $G$ and of $U^+$. A new result is that these non-negative monoids can be 
defined by explicit 
generators and relations. In this definition these monoids have non-negative $K$-structures where $K$ can be more
general than $\RR_{>0}$ (for example it can be $\ZZ$). 

In \S3 we show following \cite{L97} that the positive part $\cb_{>0}$ 
of the real flag manifold $\cb$ of $G$ (defined in \cite{L94}) has a natural positive structure. (This result
of \cite{L97} is one of the ingredients in the work of Fock and Goncharov \cite{FG} on higher Teichm\"uller 
theory.)
By passage to zones we can replace $\RR_{>0}$ by $\ZZ$ and we obtain a set $\cb(\ZZ)_{>0}$ with a positive
structure (with parametrizations by $\ZZ^\nu$). This set can be identified with $\ti\BB$ and by 3.4(d) it admits
a natural action of the monoid attached to $G$ and $\ZZ$. This can be viewed as a new symmetry property
of the canonical basis $\BB$ (which is visible only after $\BB$ is enlarged to $\ti\BB$). In \S3 we show also 
that there is a version of the positive structure on $\cb_{>0}$ where $\cb$ is replaced by $G/U^-$ (here $U^-$
is the strictly lower triangular part of $G$). 

In \S4 we study two involutions $\Ph,\Ph'$ of $\ti\BB$ and their lifting over $\RR_{>0}$ which originate in the
study \cite{L97} of the positive structure on $\cb_{>0}$.

In \S5 we give another version of the positive structure on $\cb_{>0}$ where $\cb$ is replaced by a partial 
flag manifold; this complements the results of \cite{L98}.

As shown in \cite{L18}, a Chevalley group over an algebraically closed field of any characteristic can be 
reconstructed from the non-negative monoid attached to $G$ and $\RR_{>0}$. But this did not include
a description in the same spirit of the coordinate ring of that group; this is done (conjecturally) in \S6,
where we also describe in the same spirit the Weyl modules of that group.

In \S7 we reformulate some results of Rietsch \cite{R97} on the connected components of the intersection of two 
opposed Bruhat cells in $\cb$ in terms the non-negative monoid attached to $U^+$ and $K=\ZZ$. The result 
in \S7 suggests that at least in some cases, the set of connected components of a real algebraic variety
can be described in terms of a positive structure involving $\RR(t)$.

In \S8 we are concerned with the subset $\BB_\l$ of $\BB$ attached in \cite{L90} to a dominant weight $\l$; this 
subset has the property that, when applied to a lowest weight vector in the simple $\UU$-module $V_\l$
indexed by $\l$, it gives rise to a basis of $V_\l$. (The elements in $\BB-\BB_\l$ applied to that lowest weight 
vector give the zero vector.) In \cite{L90}, the subset $\BB_\l$ was given a purely combinatorial description
in terms of the piecewise linear structure of $\BB$. In \S8 we give another combinatorial description of $\BB_\l$ 
(this time conjectural). One of new ingredients in this description is the involution $\Ph$ in \S4.

In \S9 we prove some (partly conjectural) properties of the semisimple and unipotent part of an element in 
the non-negative monoid attached to $G$ and $\RR_{>0}$.

In \S10 we define a partition of $\BB$ indexed by the Weyl group.

\subhead 0.3\endsubhead
In this paper we write 
$$U^+(K),U^+_{w_I}(K),U^-(K),U^-_{w_I}(K),G(K),G_{w_I,-w_I}(K),\cb(K)_{>0}$$ 
with $K=\RR_{>0}$ for what in \cite{L94} was denoted by 
$$U^+_{\ge0},U^+_{>0},U^-_{\ge0},U^-_{>0},G_{\ge0},G_{>0},\cb_{>0}.$$ 
We use the adjective ``positive'' (resp. ``non-negative'') for what in \cite{P10} is called
``stricly totally positive'' (resp. ``totally positive'').

In this paper, a monoid is always understood to have $1$.

\subhead 0.4\endsubhead
{\it Errata to \cite{L97}}. In Theorem 2.10 delete the first sentence. In 2.10(a),(b), 4.10, 4.11, replace $f^*$ 
by $f$.

\subhead 0.5\endsubhead
I thank Xuhua He and Konstanze Rietsch for discussions.

\head Contents \endhead
1. Positive $K$-structures.

2. The monoids $\fU(K),\fG(K)$.

3. The positive part of a flag manifold.

4. The involutions $\Ph,\Ph'$ of $\fU_{w_I}(K)$.

5. The positive part of a partial flag manifold.

6. Coordinate rings.

7. Arnold's problem.

8. The sets $\fU^\l_{w_I}(\cn),\ti\fU^\l_{w_I}(\cn)$.

9. On $G(K)$ and conjugacy classes.

10. A partition of $\BB$.

\head 1. Positive $K$-structures\endhead
\subhead 1.1\endsubhead
The non-negative monoid in 0.2 can be defined not only over $\RR_{>0}$ but over a structure $K$ in which addition, 
multiplication, division (but no substraction) are defined. In \cite{L94} two types of such $K$ were considered.

(i) There exists a field $\kk$ (necessarily of characteristic zero) such that $K\sub\kk-\{0\}$ and the addition, 
multiplication, division on $K$ are induced from the analogous operations in $\kk$. (We can take, for example, 
$\kk=\RR,K=\RR_{>0}$ or $\kk=\RR(t)$ with $t$ an indeterminate and $K=\RR(t)_{>0}$ to be the set of $f\in\kk$ of 
form $f=t^ef_0/f_1$ for some $f_0,f_1$ in $\RR[t]$ with constant term in $\RR_{>0},e\in\ZZ$.)
 
(ii) $K=\ZZ$ with a new sum $(a,b)\m\min(a,b)$ and a new product $(a,b)\m a+b$. (We write $\cz$ for $\ZZ$ with 
these new operations.) 

A third type of $K$ is:

(iii) $K=\{1\}$ with $1+1=1,1\T1=1$.
\nl
In each case $K$ is a semifield (a terminology of Berenstein, Fomin, Zelevinsky \cite{BFZ}): a set with two 
operations, $+$, $\T$, which is an abelian group with respect to $\T$, an abelian semigroup with respect to $+$ 
and in which $(a+b)c=ac+bc$ for all $a,b,c$.  In this paper, $K$ is as in (i)-(iii). (However, most definitions
and results of this paper remain valid for any semifield.) There is an obvious semifield homomorphism 
$K\to\{1\}$. We denote by $K_0$ the smallest semifield contained in $K$ which contains the unit element $(1)$ of 
$K$ with respect to $\T$. If $K$ is as in (i), we have $K_0=\QQ_{>0}$. If $K$ is as in (ii) we have $K_0=\{0\}$. 

\subhead 1.2\endsubhead
For any $m\in\ZZ_{>0}$ let $\fP_m$ be the set of rational functions $\ph=f/f'$ where $f,f'$ are nonzero
polynomials in the indeterminates $X_1,X_2,\do,X_m$ with coefficients in $\NN$. A map $K^m@>>>K^{m'}$ (with 
$m\in\ZZ_{>0},m'\in\ZZ_{>0}$) is said to be {\it admissible} if it is of the form
$(a_1,a_2,\do,a_m)\m(\ph_1(a_1,\do,a_m),\ph_2(a_1,\do,a_m),\do,\ph_{m'}(a_1,\do,a_m))$
for some $\ph_1,\ph_2,\do,\ph_{m'}\in\fP_m$. (Note that $\ph_i(a_1,\do,a_m)$ is a well defined element of $K$ if 
$(a_1,\do,a_m)\in K^m$.) In the case where $K=\ZZ$, such a map is piecewise-linear. We interpret $K^0$ to be a 
point. If $m\in\NN$, the unique map $K^m@>>>K^0$ is considered to be admissible. If $m'\in\ZZ_{>0}$, a map 
$K^0@>>>K^{m'}$ is said to be admissible if its image is a point in $K_0^{m'}$. A bijection $K^m@>\si>>K^m$ (with
$m\in\NN$) is said to be {\it bi-admissible} if it is admissible and its inverse is admissible. From the 
definitions we have the following result (with $m,m',m''$ in $\NN$).

(a) {\it A composition of admissible maps $K^m@>>>K^{m'}$, $K^{m'}@>>>K^{m''}$ is an admissible map 
$K^m@>>>K^{m''}$.}

\subhead 1.3\endsubhead
A {\it positive $K$-structure} on a set $X$ consists of a family of bijections $f_j:K^m@>\si>>X$ (with 
$m\ge0$ fixed) indexed by $j$ in a finite nonempty set $\cj$, such that $f_{j'}\i f_j:K^m@>>>K^m$ 
is bi-admissible for any $j,j'$ in $\cj$; the bijections $f_j$ are said to be the {\it charts} of the positive
structure. A {\it non-negative $K$-structure} on a set $X$ is a partition $X=\sqc_{\c\in H}X_\c$ (with 
$H$ finite) together with a positive $K$-structure 
on $X_\c$ for every $\c\in H$; the subsets $X_\c$ are said to be the {\it pieces} of $X$.

If $\{f_j:K^m@>>>X,j\in\cj\}$ (resp. $\{f'_{j'}:K^{m'}@>>>X',j'\in\cj'\}$) is a positive $K$-structure on a set 
$X$ (resp. $X'$) then $\{f_j\T f_{j'}:K^{m+m'}@>>>X\T X',(j,j')\in\cj\T\cj'\}$ is a positive $K$-structure on 
$X\T X'$; a map $\ph:X@>>>X'$ is said to be admissible  
if for some (or equivalently any) $(j,j')\in\cj\T\cj'$, 
the map $f_{j'}\i\ph f_j:K^m@>>>K^{m'}$ is admissible. If $X$ (resp. $X'$) is a set with 
non-negative 
$K$-structure with pieces $X_\c,\c\in H$ (resp. $X'_{\c'},\c'\in H'$) then $X\T X'$ has a non-negative 
$K$-structure with pieces $X_\c\T X_{\c'},(\c,\c')\in H\T H'$; a map $\ph:X@>>>X'$ is said to be admissible if 
for any $\c\in H$, $\ph(X_\c)$ is contained in $X'_{\c'}$ for a (well defined) $\c'\in H'$ and the restriction of
$\ph$ from $X_\c$ to $X'_{\c'}$ is admissible.

When $K$ is as in 1.1(iii), a set with positive $K$-structure is just a point with no structure; a set with
 non-negative $K$-structure is just a finite set with no structure.

\subhead 1.4\endsubhead
We now assume that $K=\RR(t)_{>0}$. In \cite{L94} we observed that there is a semifield homomorphism 
$\a:K@>>>\cz$ given by $t^ef_0/f_1\m e$ which connects geometrical objects over $K$ with piecewise linear 
objects involving only integers. 

Let $\{f_j:K^m@>\si>>X,j\in\cj\}$ be a positive $K$-structure on a set $X$. Assume first that $m>0$.
For any $j\in\cj$ we define an equivalence relation on $X$ in which $x,x'$ in $X$ are equivalent
if we have $f_j\i(x)=(a_1,a_2,\do,a_m)$, $f_j\i(x')=(a'_1,a'_2,\do,a'_m)$ where $a_k,a'_k$ in $K$ satisfy 
$\a(a_k)=\a(a'_k)$ for $k=1,\do,m$. If this condition holds for some $j\in\cj$ then it holds for any $j'\in\cj$. 
Indeed, setting $f_{j'}\i(x)=(b_1,b_2,\do,b_m)$, $f_{j'}\i(x')=(b'_1,b'_2,\do,b'_m)$ with $b_k,b'_k$ in $K$ we 
have 

$(b_1,b_2,\do,b_m)=(\ph_1(a_1,\do,a_m),\ph_2(a_1,\do,a_m),\do,\ph_m(a_1,\do,a_m))$,

$(b'_1,b'_2,\do,b'_m)=(\ph_1(a'_1,\do,a'_m),\ph_2(a'_1,\do,a'_m),\do,\ph_m(a'_1,\do,a'_m))$
\nl
where $\ph_k\in\fP_m$. Since $\a$ is a semifield homomorphism it follows that 
$$\align&(\a(b_1),\a(b_2),\do,\a(b_m))\\&
=(\ph_1(\a(a_1),\do,\a(a_m)),\ph_2(\a(a_1),\do,\a(a_m)),\do,\ph_m(\a(a_1),\do,
\a(a_m)))\\&=
(\ph_1(\a(a'_1),\do,\a(a'_m)),\ph_2(\a(a'_1),\do,\a(a'_m)),\do,\ph_m(\a(a'_1),\do,\a(a'_m)))\\&=
(\a(b'_1),\a(b'_2),\do,\a(b'_m)).\endalign$$
so that $\a(b_k)=\a(b'_k)$ for $k=1,\do,m$. 
Thus, our equivalence relation does not depend on the choice of $j\in\cj$. Following
\cite{L94, \S9}, the equivalence classes are called the {\it zones} of $X$.
When $m=0$, $X$ consists of one point and by definition it is a single zone.
Let $\un X$ be the set of zones of $X$. If $m>0$, for any $j\in\cj$,  we define $\un f_j:\cz^m@>>>\un X$ by 
$$(c_1,c_2,\do,c_m)\m \text{zone of $X$ containing } f_j(a_1,a_2,\do,a_m)$$
where $(a_1,a_2,\do,a_m)\in K^m$ is such that $\a(a_k)=c_k$ for $k=1,\do,m$. This map is well defined by the
definition of zones and is a bijection. Moreover, for $j,j'$ in $\cj$, the bijection
$\un f_{j'}\i\un f_j:\cz^m@>>>\cz^m$ is bi-admissible; it is induced by the bi-admissible bijection
$f_{j'}\i f_j:K^m@>>>K^m$. Thus $\{\un f_j:\cz^m@>>>\un X;j\in\cj\}$ is a positive $\cz$-structure on $\un X$.
The same is true if $m=0$ when $\un f_j$ are taken as bijections from a point to a point.

If $X,X'$ are two sets with positive $K$-structure, then $\un{X\T X'}=\un X\T\un X'$ as sets with positive 
$\cz$-structure; if in addition $\ph:X@>>>X'$ is an admissible map, then $\ph$ maps any zone of $X$ into a
zone of $X'$ and the induced map $\un\ph:\un X@>>>\un X'$ is admissible.

\subhead 1.5\endsubhead
In the setup of 1.4, let $X$ be a set with a given non-negative $K$-structure with pieces $X_\c,\c\in H$. Let
$\un X=\sqc_{\c\in H}\un X_\c$. By 1.4, each $\un X_\c$ has a positive $\cz$-structure. This defines a 
non-negative $\cz$-structure on $\un X$ with pieces $\un X_\c$, $\c\in H$.

If $X,X'$ are two sets with non-negative $K$-structure with pieces $X_\c,\c\in H$ and $X'_{\c'},\c'\in H'$, then 
$\un {X\T X'}=\un X\T\un X'$ as sets with non-negative $\cz$-structure; if in addition $\ph:X@>>>X'$ is an 
admissible map then there is a well defined admissible map $\un\ph:\un X@>>>\un X'$ such that the following holds:
if $\c\in H$ and $\c'\in H'$ is such that $\ph(X_\c)\sub X'_{\c'}$ then $\un\ph(\un{X_\c})\sub\un{X'_{\c'}}$ and 
the map $\un\ph:\un{X_\c}@>>>\un{X'_{\c'}}$ is induced (as in 1.4) by $\ph:X_\c@>>>X'_{\c'}$.

\head 2. The monoids $\fU(K),\fG(K)$\endhead
\subhead 2.1\endsubhead
Let $\kk$ be an infinite field. We will often identify an algebraic variety over $\kk$ with its set of 
$\kk$-points. We fix a split, connected, simply connected algebraic group $G$ over $\kk$ with a given $\kk$-split 
maximal torus $T$ and a pair $B^+,B^-$ of Borel subgroups with $B^+\cap B^-=T$, with unipotent radicals 
$U^+,U^-$. Let $NT$ be the normalizer of $T$ in $G$ and let $W=NT/T$ be the Weyl group.
Let $\cy$ (resp. $\cx$) be the abelian group of all homomorphisms of algebraic groups $\kk^*@>>>T$ (resp.
$T@>>>\kk^*$) with operation written as $+$. We define a pairing $\la,\ra:\cy\T\cx@>>>\ZZ$ by 
$x(y(a))=a^{\la y,x\ra}$ for all $a\in\kk^*$. The set of simple coroots (resp. simple roots) determined by 
$B^+,T$ is a $\ZZ$-basis $I$ of $\cy$ (resp. a subset $\{i^*;i\in I\}$ of $\cx$). We set $r=\sha(I)$. The 
matrix $(\la i,j^*\ra)$ indexed by $I\T I$ is the Cartan matrix of $G$. For $i\in I$, let $U_i^+$ be the simple 
root subgroup of $U^+$ defined by $i^*$ and let $U_i^-$ ($i\in I$) be the corresponding root subgroup of $U^-$.
We assume that for any $i\in I$ we are given  isomorphisms $a\m x_i(a)$, $\kk@>\si>>U_i^+$ and $a\m y_i(a)$, 
$\kk@>\si>>U_i^-$ as algebraic groups such that the assignment
$$\left(\sm1&a\\0&1\esm\right)\m x_i(a),\left(\sm b&0\\0&b\i\esm\right)\m i(b),
\left(\sm1&0\\c&1\esm\right)\m y_i(a)$$
where $a,c\in\kk,b\in\kk^*$, defines a homomorphism $SL_2(\kk)@>>>G$. Let $\ds_i$ be the image of 
$\left(\sm0&1\\-1&0\esm\right)$ under this homomorphism. We have $\ds_i\in NT$. For $a,b$ in $\kk^*$ we have
$$x_i(a)i(b)=y_i(a\i)i(ab\i)\ds_iy_i(a\i b^2).\tag a$$
Let $s_i$ be the image of $\ds_i$ in $W$. Then $W$ is a Coxeter group on the generators $\{s_i;i\in I\}$. Let 
$w\m|w|$, $W@>>>\NN$ be the standard length function. For $J\sub I$ let $W_J$ be the subgroup of $W$ generated 
by $\{i;i\in J\}$. Let $w_J$ be the element of $W_J$ such that $|w_J|$ is maximal. In particular, $w_I\in W$ is 
defined. We set $\nu=|w_I|$. Let $i\m i^!$ be the involution $I@>>>I$ such that $w_Is_iw_I=s_{i^!}$ for all 
$i\in I$. There is a unique $W$-action on $\cx$ such that for $i\in I,\l\in\cx$ we have 
$s_i(\l)=\l-\la i,\l\ra i^*$ and a unique $W$-action on $\cy$ such that for $i\in I,\z\in\cy$ we have 
$s_i(\z)=\z-\la\z,i^*\ra i$. 

We say that $G$ is {\it simply laced} if for any $i\ne j$ in $I$, $s_is_j$ has order $\le3$ in $W$.

For $w\in W$ with $m=|w|$ let $\co_w$ be the set of all sequences $(i_1,i_2,\do,i_m)$ such that 
$s_{i_1}s_{i_2}\do s_{i_m}=w$. Following Tits, for $w\in W$ with $m=|w|$, we set 
$\dw=\ds_{i_1}\ds_{i_2}\do\ds_{i_m}$ where $\ii=(i_1,i_2,\do,i_m)\in\co_w$; this is independent of the choice of
$\ii$. Let $G_w=B^-\dw B^-\sub G$. 

There is a well defined isomorphism $\io:G@>>>G$ such that $\io(x_i(a))=y_i(a)$, $\io(y_i(a))=x_i(a)$ for 
$i\in I,a\in\kk$. There is a well defined isomorphism $\Ps:G@>>>G^{opp}$ (the opposite group structure) such that
$\Ps(x_i(a))=y_i(a)$, $\Ps(y_i(a))=x_i(a)$ for $i\in I,a\in\kk$. Let $\o:G@>>>G$ be the isomorphism such that
$\o(x_i(a))=y_i(a)$, $\o(y_i(a))=x_i(a)$ for all $i\in I,a\in\kk$. For 
$i\in I$ we have $\Ps(\ds_i)=\ds_i\i$, $\o(\ds_i)=\ds_i\i$.

\mpb

Let $\cx^+=\{\l\in\cx;\la i,\l\ra\in\NN \qua\frl i\in I\}$. Let $\cc$ be the category whose objects are finite 
dimensional $\kk$-vector spaces with a given rational $G$-module structure.
For any $\l\in\cx^+$ we denote by $\L_\l$ a simple object of $\cc$ with a nonzero vector $\et_\l$
which is fixed by the $U^+$-action and satisfies $t\et_\l=\l(t)\et_\l$ for all $t\in T$. 

\subhead 2.2\endsubhead
We now assume that $\kk,K$ are as in 1.1(i). The following identities hold:
$$tx_i(a)=x_i(i^*(t)a)t,ty_i(a)=y_i(i^*(t)\i a)t\text{ for }i\in I,t\in T,a\in K;$$
$$y_i(a)x_j(c)=x_j(c)y_i(a)\text{ for $a,c\in K$ and $i\ne j$ in $I$};$$
$$y_i(a)i(b\i)x_i(c)=x_i(c')i(b')y_i(a')$$
for $i\in I$ and $a,b,c,a',b',c'$ in $K$ such that 
$$a'=a/(ac+b^2),b'=b/(ac+b^2),c'=c/(ac+b^2)$$
or equivalently $a=a'/(a'c'+b'{}^2),b=b'/(a'c'+b'{}^2),c=c'/(a'c'+b'{}^2).$

\subhead 2.3\endsubhead
Let $w\in W$ with $m=|w|$. The statements (a),(b) below are proved in \cite{L94, 2.7(a),(b)}.

(a) {\it Let $\ii=(i_1,i_2,\do,i_m)\in\co_w$. The map $\t^+_\ii:K^m@>>>U^+$ given by
$$(a_1,a_2,\do,a_m)\m x_{i_1}(a_1)x_{i_2}(a_2)\do x_{i_m}(a_m)$$ 
is injective.}

(b) {\it Let $\ii=(i_1,i_2,\do,i_m)\in\co_w,\ii'=(i'_1,i'_2,\do,i'_m)\in\co_w$. There is a well defined 
bijection $R_\ii^{\ii'}:K^m@>>>K^m$ such that for any $(a_1,a_2,\do,a_m)\in K^m$ we have
$$x_{i_1}(a_1)x_{i_2}(a_2)\do x_{i_m}(a_m)=x_{i'_1}(a'_1)x_{i'_2}(a'_2)\do x_{i'_m}(a'_m)$$
where $(a'_1,a'_2,\do,a'_m)=R_\ii^{\ii'}(a_1,a_2,\do,a_m)\in K^m$ that is, 
$\t^+_{\ii'}R_\ii^{\ii'}=\t^+_\ii:K^m@>>>K^m$.}
\nl
Now let $\ii,\ii'$ be as in (b) and let $\ii''=(i''_1,i''_2,\do,i''_m)\in\co_w$. For $(a_1,a_2,\do,a_m)\in K^m$ 
we have
$$\align&x_{i_1}(a_1)x_{i_2}(a_2)\do x_{i_m}(a_m)=x_{i'_1}(a'_1)x_{i'_2}(a'_2)\do x_{i'_m}(a'_m)\\&=
x_{i''_1}(a''_1)x_{i''_2}(a''_2)\do x_{i''_m}(a''_m)\endalign$$
where 
$$\align&(a'_1,a'_2,\do,a'_m)=R_\ii^{\ii'}(a_1,a_2,\do,a_m),(a''_1,a''_2,\do,a''_m)=R_{\ii'}^{\ii''}(a'_1,a'_2,\do,a'_m),\\&
(a''_1,a''_2,\do,a''_m)=R_{\ii}^{\ii''}(a_1,a_2,\do,a_m).\endalign$$

It follows that 
$$R_\ii^{\ii''}(a_1,a_2,\do,a_m)=R_{\ii'}^{\ii''}(a'_1,a'_2,\do,a'_m)=
R_{\ii'}^{\ii''}R_\ii^{\ii'}(a_1,a_2,\do,a_m).$$
Hence 

(c) $R_\ii^{\ii''}=R_{\ii'}^{\ii''}R_\ii^{\ii'}$.
\nl
From the proof in \cite{L94, 2.6, 2.7} one can also extract the following result.

(d) {\it The bijection $R_\ii^{\ii'}:K^m@>>>K^m$ is bi-admissible.}
\nl
By a lemma of Iwahori \cite{Iw}, 
$\ii,\ii'$ can be connected by a finite sequence of braid moves. Using this, (c)
and 
1.2(a), we see that to prove (d) we can assume that $\ii,\ii'$ are connected by a braid move that is, $\ii'$ is 
obtained from $\ii$ by replacing $e$ consecutive entries
$i,j,i,\do$ of $\ii$ by the $e$ entries $j,i,j,\do$ where $i\ne j$ in $I$ are such that $s_is_j$ has order
$e$ in $W$. If $e$ is $2$ or $3$, the desired result follows from \cite{L94, 2.5(a),(b)}. (Hence (d) holds 
whenever $G$ is simply laced.) Assume now that $e>3$, that is $e=2e'$ with $e'\in\{2,3\}$. We can assume 
$G=G_{i,j}$, the subgroup of $G$ generated by $U_i^+,U_j^+,U_i^-,U_j^-$, so that
that $I=\{i,j\}$ and $\ii=(i,j,i,\do),\ii'=(j,i,j,\do)$ (both sequences consist of $2e'$ terms).
We can find a group $\dot G$ group like $G$ in 2.1 but of type $A_3$ (if $e'=2$) or $D_4$ (if $e'=3$) and
an automorphism of order $e'$ of $\dot G$ whose fixed point set is $G$. We have $\la i,j^*\ra\in\{-1,-e'\}$. Let 
$\dW$ be the Weyl group of $\dot G$. 
The simple reflections of $\dW$ are denoted $\s_1,\s_2,\do,\s_{e'},\s_0$ where
$\s_a,\s_b$ commute for $a,b\in\{1,2,\do,e'\}$ and $\s_a\s_0$ has order $3$ for $a\in\{1,2,\do,e'\}$. 
(The indexing set for the simple reflections of $\dW$ is $\dI=\{1,2,\do,e',0\}$.) We can identify $W$ with a 
subgroup of $\dW$ in such a way that either

(i) $s_i$ becomes $\s_1\s_2\do\s_{e'}$, $s_j$ becomes $\s_0$ (in this case $\la i,j^*\ra=-e'$), or

(ii) $s_i$ becomes $\s_0$, $s_j$ becomes $\s_1\s_2\do\s_{e'}$ (in this case $\la i,j^*\ra=-1$). 
\nl
We consider two reduced expressions for the longest element (of length $e'(e'+1)$) of $\dW$ namely
$(1,2,\do,e',0,1,2,\do,e',0,\do)$, where $1,2,\do,e',0$ is repeated $e'$ times, and
$(0,1,2,\do,e',0,1,2,\do,e',\do)$, where $0,1,2,\do,e'$ is repeated $e'$ times.
In case (i), the first of these reduced expressions is denoted by $\dot\ii$ and the second by $\dot\ii'$.
In case (ii), the first of these reduced expressions is denoted by $\dot\ii'$ and the second by $\dot\ii$.
In case (i) we define $\l:K^{2e'}@>>>K^{e'(e'+1)}$ by
$$(a_1,a_2,\do,a_{2e'})\m(a_1,a_1,\do,a_1,a_2,a_3,a_3,\do,a_3,a_4,\do,a_{2e'-1},\do a_{2e'-1},a_{2e'})$$
(with $a_{2t+1}$ repeated $e'$ times) and $\mu:K^{e'(e'+1)}@>>>K^{2e'}$ by
$$(a_1,a_2,\do,a_{e'(e'+1)})\m(a_1,a_2,a_{e'+2},a_{e'+3},a_{2e'+3},a_{2e'+4},\do).$$
In case (ii) we define $\l:K^{2e'}@>>>K^{e'(e'+1)}$ by
$$(a_1,a_2,\do,a_{2e'})\m(a_1,a_2,\do,a_2,a_3,a_4,\do,a_4,\do,a_{2e'-1},a_{2e'},\do,a_{2e'})$$
(with $a_{2t}$ repeated $e'$ times) and $\mu:K^{e'(e'+1)}@>>>K^{2e'}$ by
$$(a_1,a_2,\do,a_{e'(e'+1)})\m(a_1,a_{e'+1},a_{e'+2},a_{2e'+2},a_{2e'+3},a_{3e'+3},\do).$$
We consider the composition 
$$K^{2e'}@>\l>>K^{e'(e'+1)}@>R_{\dot\ii}^{\dot\ii'}>>K^{e'(e'+1)}@>\mu>>K^{2e'}$$
where the middle map is defined in terms of $\dot G$. By the earlier part of the proof, the middle map is 
admissible (since $\dot G$ is simply laced); it follows that the composition is admissible. By the argument in 
the proof of \cite{L94, 2.6, 2.7} this composition is equal to $R_\ii^{\ii'}:K^{2e'}@>>>K^{2e'}$. Thus, 
$R_\ii^{\ii'}:K^{2e'}@>>>K^{2e'}$ is admissible. The inverse of $R_\ii^{\ii'}:K^{2e'}@>>>K^{2e'}$ is 
$R_{\ii'}^\ii:K^{2e'}@>>>K^{2e'}$ which is admissible. Thus, $R_\ii^{\ii'}:K^{2e'}@>>>K^{2e'}$ is bi-admissible 
and (d) holds.

\mpb

Another proof of (d) was given later by Berenstein and Zelevinsky, see Theorem 1.6 and 3.1 of \cite{BZ}.

\mpb

Using the automorphism $\o$, from (a),(b) we deduce:

(e) {\it Let $\ii=(i_1,i_2,\do,i_m)\in\co_w$. The map $\t^-_\ii:K^m@>>>U^-$ given by
$$(a_1,a_2,\do,a_m)\m y_{i_1}(a_1)y_{i_2}(a_2)\do y_{i_m}(a_m)$$ is injective.}

(f) {\it Let $\ii=(i_1,i_2,\do,i_m)\in\co_w,\ii'=(i'_1,i'_2,\do,i'_m)\in\co_w$. Then the bijection 
$R_\ii^{\ii'}:K^m@>>>K^m$ in (b) has the following property: for any $(a_1,a_2,\do,a_m)\in K^m$ we have
$$y_{i_1}(a_1)y_{i_2}(a_2)\do y_{i_m}(a_m)=y_{i'_1}(a'_1)y_{i'_2}(a'_2)\do y_{i'_m}(a'_m)$$
 where $(a'_1,a'_2,\do,a'_m)=R_\ii^{\ii'}(a_1,a_2,\do,a_m)\in K^m$ that is, 
$\t^-_{\ii'}R_\ii^{\ii'}=\t^-_\ii:K^m@>>>K^m$.}

\subhead 2.4\endsubhead
Let $i\ne j$ be two elemens of $I$ and let $e$ be the order of $s_is_j$ in $W$. Let $G_{i,j}$ be as in 2.3. Then 
the bijection

(a) $R(i,j):=R_{i,j,i,\do}^{j,i,j,\do}:K^e@>>>K^e$ 
\nl
(where $i,j,i,\do$ has $e$ terms and $j,i,j\do$ has $e$ terms) is well defined in terms of $G_{i,j}$. It is 
bi-admissible by 2.3(d). If $e=2$, then $R(i,j)$ is given by 

(b) $(a,b)\m(b,a)$. 
\nl
If $e=3$, then $R(i,j)$ is given by

(c) $(a,b,c)\m(a',b',c')$ where $a'=bc/(a+c),b'=a+c,c'=ab/(a+c)$ or equivalently
$a=b'c'/(a'+c'),b=a'+c',c=a'b'/(a'+c')$.
\nl
A formula like (c) appeared in \cite{L90} (for $K=\cz$) in connection with the problem of
parametrizing canonical bases.

When $e=2e'$ is $4$ or $6$ then, with the notation of 2.3, $R(i,j)$ is a product over a set of braid moves
connecting $(1,2,\do,e',0,1,2,\do,e',0,\do)$ with $(0,1,2,\do,e',0,1,2,\do,e',\do)$ in a Weyl group of type 
$A_3$ or $D_4$ of bijections of type (b) and (c) involving two or three coordinates. For example, if $e=4$ and 
$\la i,j^*\ra=-2$, then $R(i,j)$ is given by

(d) $(a,b,c,d)\m(a',b',c',d')$ where $a'=bc^2d/E,b'=E/A,c'=A^2/E,d'=abc/A$, $A=ab+ad+cd$, $E=a^2b+a^2d+c^2d+2acd$,
\nl
or equivalently

(e) $d=a'b'{}^2c'/E',c=E'/A',b=A'{}^2/E',a=b'c'd'/A'$, $A'=c'd'+a'd'+a'b'$, 
$E'=d'{}^2c'+a'd'{}^2+a'b'{}^2+2a'b'd'$.
\nl
Note that in (d),(e) we have 

(f) $abc=b'c'd',bc^2d=a'b'{}^2c',a+c=b'+d',b+d=a'+c'$.
\nl
Conversely, if $a,b,c,d,a',b',c',d'$ in $K$ satisfy (f), then either the equalities in (d) hold or we have 
$(a',b',c',d')=(d,c,b,a)$.

A formula like (d) appeared (for $K=\cz$) in \cite{L92} in connection with the problem of
parametrizing canonical bases.

\subhead 2.5\endsubhead
Let $U^+(K)$ be the submonoid of $U^+$ generated by $\{x_i(a);i\in I,a\in K\}$. Let $U^-(K)$ be the submonoid of 
$U^-$ generated by $\{y_i(a);i\in I,a\in K\}$. Let $T(K)$ be the submonoid of $T$ generated by 
$\{i(a);i\in I,a\in K\}$; it is a subgroup of $T$ since $i(a)\i=i(a\i)$ for $i\in I,a\in K$. Let $G(K)$ be the 
submonoid of $G$ generated by $\{x_i(a),y_i(a),i(a);i\in I,a\in K\}$. (These definitions appeared in 
\cite{L94, 2.2}.) We have $\Ps(U^+(K))=U^-(K)$, $\Ps(U^-(K))=U^+(K)$, $\Ps(T(K))=T(K)$, $\Ps(G(K))=G(K)$ and
$\o(U^+(K))=U^-(K)$, $\o(U^-(K))=U^+(K)$, $\o(T(K))=T(K)$, $\o(G(K))=G(K)$.
The following is proved in \cite{L94, 2.3}.

(a) {\it If $u^+\in U^+(K),t\in T(K),u^-\in U^-(K)$ then $u^+tu^-\in G(K)$ and $u^-tu^+\in G(K)$. Any $g\in G(K)$
can be written uniquely in the form $g=u^+tu^-$ with $u^+\in U^+(K),t\in T(K),u^-\in U^-(K)$. Any $g\in G(K)$ can
be written uniquely in the form $g=u_1^-t_1u_1^+$ with $u_1^+\in U^+(K),t_1\in T(K),u_1^-\in U^-(K)$.}
\nl
Let $w\in W$, $m=|w|$. Let $U^+_w(K)=\t^+_\ii(K^m)$ where $\ii\in\co_w$; by 2.3(b), this is independent of the 
choice of $\ii$ and $\t^+_\ii$ defines a bijection $K^m@>\si>>U^+_w(K)$. Let $U^-_w(K)=\t^-_\ii(K^m)$ where 
$\ii\in\co_w$; by 2.3(f), this is independent of the choice of $\ii$ and $\t^-_\ii$ defines a bijection 
$K^m@>\si>>U^-_w(K)$. Note that $U^-_w(K)=\o(U^+_w(K))$. 

Following \cite{L94, 2.11}, for $w,w'$ in $W$ we set 

(b) $G_{w,-w'}(K)=U^-_{w'}(K)T(K)U^+_w(K)=U^+_w(K)T(K)U^-_{w'}(K)$. 
\nl
(The last equality follows from repeated applications of the identities in 2.2.) 

\mpb

For future reference we state:

(c) $U^+_w(K)\sub B^-\dw B^-$.
\nl
(See the proof of \cite{L94, 2.7}).

\mpb

From \cite{L94, 2.7(a), 2.8} we see that 

(d) $U^+(K)=\sqc_{w\in W}U^+_w(K)$.
\nl
Applying $\o$ to (d) we obtain

(e) $U^-(K)=\sqc_{w\in W}U^-_w(K)$.

\subhead 2.6\endsubhead
From 2.5(a),(d),(e) we deduce

(a) $G(K)=\sqc_{w,w'\text{ in }W}G_{w,-w'}(K)$.

\subhead 2.7\endsubhead
Let $w\in W, m=|w|$. From 2.3(d) we deduce:

(a) {\it The bijections $K^m@>>>U^+_w(K)$, 
$$(a_1,a_2,\do,a_m)=x_{i_1}(a_1)x_{i_2}(a_2)\do x_{i_m}(a_m)$$
(for various $\ii=(i_1,i_2,\do,i_m)\in\co_w$) form a positive $K$-structure on $U^+_w(K)$. The bijections 
$K^m@>>>U^-_w(K)$, $$(a_1,a_2,\do,a_m)=y_{i_1}(a_1)y_{i_2}(a_2)\do y_{i_m}(a_m)$$ (for various 
$\ii=(i_1,i_2,\do,i_m)\in\co_w$) form a positive $K$-structure on $U^-_w(K)$. When $m=0$ these bijections
are interpreted as the obvious bijections from a point to a point.}
\nl
These positive $K$-structures can be viewed as a non-negative $K$-structure on $U^+(K)$ (with pieces
$U^+_w(K),w\in W$) and a non-negative $K$-structure on $U^-(K)$ (with pieces $U^-_w(K),w\in W$). 

\mpb

In the remainder of this subsection we assume that $\kk=\RR(t),K=\RR(t)_{>0}$. For $w\in W$ we set
$U^+_w(\cz)=\un{U^+_w(K)}$, a set with a positive $\cz$-structure, see 1.4. We set
$U^+(\cz)=\un{U^+(K)}=\sqc_{w\in W}U^+_w(\cz)$, a set with a non-negative $\cz$-structure, see 1.5.
In 2.12 we will see that $U^+(\cz)$ is naturally a monoid.

\mpb

Let $\si_2$ be the equivalence relation on $U^+_{w_I}(\cz)$ generated by the
relation for which $x,x'$ in $U^+_{w_I}(\cz)$ are related if for some $\ii\in\co_{w_I}$,
$x,x'$ correspond under the bijection $\cz^\nu@>>>U^+_{w_I}(\cz)$ indexed by $\ii$ to
sequences $(n_1,n_2,\do,n_\nu)$, $(n'_1,n'_2,\do,n'_\nu)$ in $\ZZ^\nu$ such that
$n_s=n'_s\mod2$ for $s=1,2,\do,\nu$. Let $U^+_{w_I}(\ZZ)/2$ be the set of equivalence
classes.

\subhead 2.8\endsubhead
We assume that $K,\kk$ are as in 1.1(i).
We now explain in more detail the argument of \cite{L94, 2.11} (based on the identities in 2.2) which was used 
to prove the last equality in 2.5(b).

In addition to the set $I$ we consider the set $-I=\{-i;i\in I\}\sub\cy$ and another set $\uI=\{\ui;i\in I\}$ in 
obvious bijection with $I$. For $w,w'$ in $W$ with $m=|w|,m'=|w'|$ let $M=m+m'+r$ and let $\co_{w,-w'}$ be the set of sequences 
$(h_1,h_2,\do,h_M)$ in $I\sqc(-I)\sqc\uI$ such that the subsequence consisting of symbols in $I$ is in $\co_w$, 
the subsequence consisting of symbols in $-I$ is of the form $(-i_1,-i_2,\do,-i_{m'})$ with 
$(i_1,i_2,\do,i_{m'})\in\co_{w'}$
in $\co_{-w'}$ and the subsequence consisting of symbols in $\uI$ contains each symbol $\ui$ 
(with $i\in I$) exactly once.

Let $\G_{w,-w'}$ be the set of all pairs $(\hh,\aa)\in\co_{w,-w'}\T K^M$. We regard $\G_{w,-w'}$ as 
the set of vertices of a graph in which 

$(\hh,\aa)=((h_1,h_2,\do,h_M),(a_1,a_2,\do,a_M))\in\G_{w,-w'}$, 

$(\hh',\aa')=((h'_1,h'_2,\do,h'_M),(a'_1,a'_2,\do,a'_M))\in\G_{w,-w'}$ 

are joined if one of (i)-(iv) below holds:

(i) for some $t$ and some $\e\in\{1,-1\}$ we have $(h_t,h_{t+1},\do,h_{t+m-1})=
(\e i,\e j,\e i,\do)$ ($m$ terms), $$(h'_t,h'_{t+1},\do,h'_{t+m-1})=(\e j,\e i,\e i,\do)$$ ($m$ terms)
where $i\ne j$, $m$ is the order of $s_is_j$, 
$$(a'_t,a'_{t+1},\do,a'_{t+m-1})=R(i,j)(a_t,a_{t+1},\do,a_{t+m-1})$$ and $h'_s=h_s,a'_s=a_s$ 
for $s\n\{t,t+1,\do,t+m-1\}$.

(ii) for some $t$ we have $(h_t,h_{t+1})=(\ui,\uj)$, $(h'_t,h'_{t+1})=(\uj,\ui)$,
$(a'_t,a'_{t+1})=(a_{t+1},a_t)$ and $h'_s=h_s,a'_s=a_s$ for $s\n\{t,t+1\}$;

(iii) for some $t$ and some $\e\in\{1,-1\}$ we have $(h_t,h_{t+1})=(\uj,\e i)$,
$(h'_t,h'_{t+1})=(\e i,\uj)$, $(a'_t,a'_{t+1})=(a_t^{\e\la j,i^*\ra)}a_{t+1},a_t)$ 
where $i,j\in I$ and $h'_s=h_s,a'_s=a_s$ for $s\n\{t,t+1\}$;

(iv) for some $t$ we have $(h_t,h_{t+1},h_{t+2})=(i,\ui,-i)$, $(h'_t,h'_{t+1},h'_{t_2})=(-i,\ui,i)$,
where $i\in I$,
$$\align&(a'_t,a'_{t+1},a'_{t+2})\\&=(a_{t+2}/(a_ta_{t+2}+a_{t+1}^2),(a_ta_{t+2}+a_{t+1}^2)/a_{t+1},
a_t/(a_ta_{t+2}+a_{t+1}^2))\endalign$$
and $h'_s=h_s,a'_s=a_s$ for $s\n\{t,t+1,t+2\}$.
\nl
We regard $\co_{w,-w'}$ as the set of vertices of a graph in which $\hh,\hh'$ are joined if
for some $\aa,\aa'$ in $K^M$, $(\hh,\aa),(\hh',\aa')$ are joined in $\G_{w,-w'}$. 
The map $\G_{w,-w'}@>>>\co_{w,-w'}$, $(\hh,\aa)\m\hh$, respects the graph structures. 

This graph $\co_{w,-w'}$ is connected: if $\hh\in\co_{w,-w'}$ we can join $\hh$ (using edges in the graph
which are images of edges of type (ii),(iii),(iv)) with an $\hh'$ in which the first $m$ terms are in $I$ and the
last $m'$ terms are in $-I$  and then we
note that, by Iwahori's lemma, any two such $\hh'$ can be joined by a sequence of edges 
which are images of edges of type (i),(ii).

We define a map $\g_{w,-w'}:\G_{w,-w'}@>>>G(K)$ by

(a) $(\hh,\aa)=((h_1,h_2,\do,h_M),(a_1,a_2,\do,a_M))\m h_1^{a_1}h_2^{a_2}\do h_M^{a_M}$ 
\nl
where $h_s^{a_s}$ is $x_i(a_s)$ if $h_s=i$, is $y_i(a_s)$ if $h_s=-i$ and is $i(a_s)$ if $h_s=\ui$.

For any $\hh=(h_1,h_2,\do,h_M)\in\co_{w,-w'}$ we define $\ps_\hh:K^M@>>>G(K)$ by 
$\ps_\hh(\aa)=\g_{w,-w'}(\hh,\aa)$.

From the definitions and from 2.3(d) we see that if $\hh,\hh'$ are joined in the graph
$\co_{w,-w'}$ then $\ps_{\hh'}=\ps_\hh\s$ where $\s:K^M@>>>K^M$ is a bi-admissible bijection
and that the image of $\ps_\hh$ is equal to the image of $\ps_{\hh'}$.
(For the last statement we use the formulas in 2.2 and 2.3(b),2.3(f),2.4(a).)
Using this and the connectedness of the graph $\co_{w,-w'}$ we see that 
for any $\hh,\hh'$ in $\co_{w,-w'}$ we have $\ps_{\hh'}=\ps_\hh\ti\s$ where $\ti\s:K^M@>>>K^M$ is a
bi-admissible bijection and that the image of $\ps_\hh$ is independent of $\hh$.
Now, if $\hh\in\co_{w,-w'}$ has the first $m$ terms in $I$ and the last $m'$ terms in $-I$, then
$\ps_\hh$ defines a bijection of $K^M$ onto $U^+_w(K)T(K)U^-_{w'}(K)$ (we use 2.5(a) and 2.7(a)). 
We see that for any $\hh$, $\ps_\hh$ is a bijection of $K^M$ onto $U^+_w(K)T(K)U^-_{w'}(K)$.
If $\hh\in\co_{w,-w'}$ has the first $m'$ terms in $-I$ and the last $m$ terms in $I$, then
the image of $\ps_\hh$ is $U^-_{w'}(K)T(K)U^+_w(K)$. We see that
the last equality in 2.5(b) holds and the definition of $G_{w,-w'}$ in 2.5(b) is justified.
From the arguments above, 

(b) {\it For $\hh\in\co_{w,-w'}$, $\ps_\hh:K^M@>>>G_{w,-w'}$ is a bijection};

(c) {\it the bijections $\ps_\hh:K^M@>\si>>G_{w,-w'}(K)$ (for various $\hh\in\co_{w,-w'}$) define
a positive $K$-structure on $G_{w,-w'}(K)$;}

(d) {\it these positive $K$-structures can be viewed as a non-negative $K$-structure on $G(K)$ (with pieces
$G_{w,-w'}(K),(w,w')\in W\T W$).}
\nl
In the remainder of this subsection we assume that $\kk=\RR(t),K=\RR(t)_{>0}$. For $w,w'$ in $W$ we set
$G_{w,-w'}(\cz)=\un{G_{w,-w'}(K)}$, a set with a positive $\cz$-structure, see 1.4. We set
$G(\cz)=\un{G(K)}=\sqc_{w,w'\text{ in }W}G_{w,-w'}(\cz)$, a set with a non-negative $\cz$-structure, see 1.5.
In 2.12 we will see that $G(\cz)$ is naturally a monoid.

\subhead 2.9\endsubhead
Assume that $K$ is as in 1.1(i)-(iii).
Let $\fU(K)$ be the monoid with generators the symbols $i^a$ with $i\in I$, $a\in K$ and with relations 

(i) $i^ai^b=i^{a+b}$ for $i\in I$, $a,b$ in $K$;

(ii) $i^{a_1}j^{a_2}i^{a_3}\do=j^{a'_1}i^{a'_2}j^{a'_3}\do$ (both products have $m$ factors)
where $i\ne j$ in $I$, $m$ is the order of $s_is_j$ and $(a_1,a_2,\do,a_m)\in K^m$,
$(a'_1,a'_2,\do,a'_m)\in K^m$ are such that $(a'_1,a'_2,\do,a'_m)=R(i,j)(a_1,a_2,\do,a_m)$.
\nl
 In the case where $K=\cz$ and $G$ is simply laced, the definition of $\fU(K)$ appeared in \cite{L94, 9.11}.
The definition of $\fU(K)$ is reminiscent of the definition of the Coxeter group attached to the
Cartan matrix $(\la i,j^*\ra)$.

For any $w\in W$ with $m=|w|$ and any $\ii=(i_1,i_2,\do,i_m)\in\co_w$ 
let $\fU_w(K)$ be image of the map $e_\ii:K^m@>>>\fU(K)$ given by
$(a_1,a_2,\do,a_m)\m i_1^{a_1}i_2^{a_2}\do i_m^{a_m}$; this image is independent of the choice of $\ii$.
Indeed if $\ii,\ii'$ are in $\co_w$ then $\ii,\ii'$ are connected by a sequence of braid moves
(Iwahori's lemma) so it is enough to show that $e_\ii(K^m)=e_{\ii'}(K^m)$ when $\ii,\ii'$ are connected
by a single braid move. But in this case the desired equality follows from (ii).
Let $\fU'=\cup_{w\in W}\fU_w(K)\sub\fU(K)$. We show that

(a) $i^a\fU'\sub\fU'$
\nl
for any $i\in I,a\in K$. It is enough to show that for $w\in W$ the following holds:

(b) if $|s_iw|>|w|$ then $i^a\fU_w(K)\sub\fU_{s_iw}(K)$;

(c) if $|s_iw|<|w|$ then $i^a\fU_w(K)\sub\fU_w(K)$.
\nl
Now (b) is clear from the definition. In the setup of (c) we have $\fU_w(K)=\cup_{b\in K}i^b\fU_{s_iw}(K)$
hence $i^a\fU_w(K)\sub\cup_{b\in K}i^{a+b}\fU_{s_iw}(K)\sub\fU_w(K)$ and (c) is proved. Thus, (a) holds.
From (a) we deduce that $\fU'=\fU(K)$ that is, $\fU(K)=\cup_{w\in W}\fU_w(K)$. (This proof is almost a copy
of that in \cite{L94, 2.8}.)

For each $w\in W$ we choose $\ii_w\in\co_w$. By the arguments above, the map

$\z'=\sqc_{w\in W}e_{\ii_w}:\sqc_{w\in W}K^{|w|}@>>>\fU(K)$
\nl
 is surjective.

Assume now that $K$ is as in 1.1(i). 
There is a well defined homomorphism of monoids $\z:\fU(K)@>>>U^+(K)$ such that $i^a\m x_i(a)$
for any $i\in I,a\in K$. (We use 2.3(b),2.4(a).) 
The composition $\z\z':\sqc_{w\in W}K^{|w|}@>>>U^+(K)$ is a bijection (we use 2.3(a) and 2.5(d)) hence
$\z'$ is injective. It follows that $\z'$ is bijective hence $\z$ is bijective. Thus the following holds.

(d) {\it For $K$ as in 1.1(i), the homomorphism of monoids $\z:\fU(K)@>>>U^+(K)$ is an isomorphism.}
\nl
Similarly,

(e) {\it For $K$ as in 1.1(i), there is a unique isomorphism of monoids $\fU(K)@>\si>>U^-(K)$ such that
$i^a\m y_i(a)$ for any $i\in I,a\in K$.}

\mpb

Let $K$ be as in 1.1(i)-(iii). We show:

(f) {\it there is a unique isomorphism of monoids $\ti\Ps:\fU(K)@>>>\fU(K)^{opp}$ (the 
opposed monoid) such that $i^a\m i^a$ for any $i\in I,a\in K$.}
\nl
Assume first that $K,\kk$ are as in 1.1(i). Then there is a unique isomorphism of monoids 
$\ti\Ps:\fU(K)@>>>\fU(K)^{opp}$ such that the following diagram is commutative:
$$\CD
\fU(K)@>\ti\Ps>>\fU(K)\\
@VVV             @VVV\\
U^+(K)@>\Ps>>  U^-(K)\endCD$$
Here the vertical maps are as in (d),(e) and $\Ps$ is as in 2.1. Then (f) follows in our case.
The case where $K=\cz$ is obtained from the case where $K=\RR(t)_{>0}$ by passage to zones. The case
where $K=\{1\}$ is immediate.

\subhead 2.10\endsubhead
Assume that $K$ is as in 1.1(i)-(iii). Let $\fG(K)$ be the monoid with generators the symbols
$i^a,(-i)^a,\ui^a$ with $i\in I$, $a\in K$ and with relations (i)-(vii) below. 

(i) $(e i)^a(\e i)^b=(\e i)^{a+b}$ for $i\in I$, $\e=\pm1$, $a,b$ in $K$;

(ii) $(\e i)^{a_1}(\e j)^{a_2}(\e i)^{a_3}\do=(\e j)^{a'_1}(\e i)^{a'_2}(\e j)^{a'_3}\do$ (both products have $m$ 
factors) where $e=\pm1$, $i\ne j$ in $I$, $m$ is the order of $s_is_j$ and $(a_1,a_2,\do,a_m)\in K^m$, 
$(a'_1,a'_2,\do,a'_m)\in K^m$ are such that $(a'_1,a'_2,\do,a'_m)=R(i,j)(a_1,a_2,\do,a_m)$;

(iii) $i^a\ui^b(-i)^c=(-i)^{c/(ac+b^2)}\ui^{(ac+b^2)/b}i^{a/(ac+b^2)}$ for $i\in I$, $a,b,c$ in $K$;

(iv) $\ui^a\ui^b=\ui^{ab}$, $\ui^{(1)}=1$ for $i\in I$, $a,b$ in $K$;

(v) $\ui^a\uj^b=\uj^b\ui^a$ for $i,j$ in $I$, $a,b$ in $K$;

(vi) $\uj^a(\e i)^b=(\e i)^{a^{\e\la j,i^*\ra}b}\uj^a$ for $i,j$ in $I$, $\e=\pm1$, $a,b$ in $K$;

(vii) $(\e i)^a(-\e j)^b=(-\e j)^b(\e i)^a$ for $i\ne j$ in $I$, $\e=\pm1$, $a,b$ in $K$.
\nl
Let $w,w'$ in $W$ be such that $|w|+|w'|+r=M$. Let $\g_{w,-w'}:\G_{w,-w'}@>>>G(K)$ be as in 2.8.
We define a map $\d:\G_{w,-w'}@>>>\fG(K)$ by

(a) $(\hh,\aa)=((h_1,h_2,\do,h_M),(a_1,a_2,\do,a_M))\m h_1^{a_1}h_2^{a_2}\do h_M^{a_M}$.
\nl
Note that $\b\d=\g_{w,-w'}$.

For any $\hh=(h_1,h_2,\do,h_M)\in\co_{w,-w'}$ we define $\th_\hh:K^M@>>>\fG(K)$ by 
$\th_\hh(\aa)=\d(\hh,\aa)$. As in the proof in 2.8 we see that the image of $\th_\hh$ is independent
of the choice of $\hh$; we denote it by $\fG_{w,-w'}(K)$.

Let $\fG'=\cup_{w\in W,w'\in W}\fG_{w,-w'}(K)\sub\fG(K)$. We show

(b) $i^a\fG'\sub\fG'$
\nl
for any $i\in I,a\in K$. It is enough to show that for $w\in W$ the following holds:

(c) if $|s_iw|>|w|$ then $i^a\fG_{w,-w'}(K)\sub \fG_{s_iw,-w'}(K)$;
if $|s_iw|<|w|$ then $i^a\fG_{w,-w'}(K)\sub\fG_{w,-w'}(K)$.
\nl
The proof is entirely similar to that of 2.9(b),(c). Similarly we have

(d) $(-i)^a\fG'\sub\fG'$
\nl
for any $i\in I,a\in K$. From the definitions we have 

(e) $\ui^a\fG'\sub\fG'$

for any $i\in I,a\in K$. From (a),(d),(e) we see that $\fG'=\fG(K)$.

For each $w,w'$ in $W$ we choose $\hh_{w,-w'}\in\co_{w,-w'}$. By the arguments above,
$\b':=\sqc_{w,w'\text{ in }W}\th_{\hh_{w,-w'}}:\sqc_{w,w'\text{ in }W}K^{|w|+|w'|+r}@>>>\fG(K)$ 
is surjective. 

Assume now that $K$ is as in 1.1(i). 
There is a well defined homomorphism of monoids $\b:\fG(K)@>>>G(K)$ such that $i^a\m x_i(a)$,
$(-i)^a\m y_i(a)$, $\ui)^a\m i(a)$ for any $i\in I,a\in K$. (We use 2.2, 2.3(b), 2.3(f), 2.4(a).) 
The composition $\b\b':\sqc_{w,w'\text{ in }W}K^{|w|+|w'|+r}@>>>G(K)$ is a bijection 
(we use 2.6(a) and 2.8(b)) hence
$\b'$ is injective. It follows that $\b'$ is bijective hence $\b$ is bijective. Thus the following holds.

(f) {\it If $K$ is as in 1.1(i), the homomorphism of monoids $\b:\fG(K)@>>>G(K)$ is an isomorphism.}

\subhead 2.11\endsubhead
In this subsection $K$ is as in 1.1(i)-(iii). As in 2.9, 2.10 we have:

(a) $\fU(K)=\cup_{w\in W}\fU_w(K)$, $\fG(K)=\cup_{w,w'\text{ in }W}\fG_{w,-w'}(K)$.
\nl
Assume now that $K=\{1\}$. In this case 
$\fU(\{1\})$ is the monoid with generators the symbols $i^1=i$  with $i\in I$ and with relations 
$ii=i$ for $i\in I$ and $iji\do=jij\do$ (both products have $m$ factors) where $i\ne j$ in $I$ and
$m$ is the order of $s_is_j$. If $w\in W$, $|w|=m$, then $\fU_w(\{1\})$ consists of a single element given 
by $i_1i_2\do i_m$  where $(i_1,i_2,\do,i_m)$ is any element of $\co_w$. We have a bijection
$W@>\si>>\fU(\{1\})$ which takes $w$ to the unique element of $\fU_w(\{1\})$.

Now $\fG(\{1\})$ is the monoid with generators the symbols $i,-i$  with $i\in I$ and with relations 

$(e i)(\e i)=(\e i)$ for $i\in I$, $\e=\pm1$;

$(\e i)(\e j)(\e i)\do=(\e j)(\e i)(\e j)\do$ (both products have $m$ factors) where $i\ne j$ in $I$,
$\e=\pm1$  and $m$ is the order of $s_is_j$;

$i(-j)=(-j)i$ for $i,j\in I$.
\nl
It follows that we have an obvious isomorphism of monoids $\fU(\{1\})\T\fU(\{1\})@>\si>>\fG(\{1\})$.
For $w,w'$ in $W$, $|w|=m,|w'|=m'$,  $\fG_{w,-w'}(\{1\})$ consists of a single element given by 
$i_1i_2\do i_m(-i'_1)(-i'_2)\do(-i'_{m'})$ where $(i_1,i_2,\do,i_m)$ is any element of $\co_w$ and
$(i'_1,i'_2,\do,i'_{m'})$ is any element of $\co_{w'}$. We have a bijection $W\T W@>\si>>\fG(\{1\})$ which 
takes $(w,w')$ to the unique element of $\fG_{w,-w'}(\{1\})$. 

\mpb

For any $K$ as in 1.1(i)-(iii) we have a homomorphism of monoids $h:\fU(K)@>>>\fU(\{1\})$ given by $i^a\m i$ and 
a homomorphism of monoids $h':\fG(K)@>>>\fG(\{1\})$ given by $i^a\m i, (-i)^a\m-i$, $\ui^a\m1$.
From the definitions we have $h(\fU_w(K))=\fU_w(\{1\})$ for $w\in W$ and
$h'(\fG_{w,-w'}(K))=\fG_{w,-w'}(\{1\})$ for $w,w'$ in $W$. Since $\fU_w(\{1\})$ are disjoint for
various $w\in W$, it follows that $\fU_w(K)$ are disjoint for various $w\in W$.
 Since $\fG_{w,-w'}(\{1\})$ are disjoint for various $w,w'$ in $W$, it follows that 
$\fG_{w,-w'}(K)$ are disjoint for various $w,w'$ in $W$. Combining this with (a) we obtain:

(b) $\fU(K)=\sqc_{w\in W}\fU_w(K)$, $\fG(K)=\sqc_{w,w'\text{ in }W}\fG_{w,-w'}(K)$.
\nl
(This is already known for $K$ as in 1.1(i).) Since $h,h'$ above are monoid homomorphisms we see that:

(c) {\it for $w,w'$ in $W$, the multiplication in $\fU(K)$ maps $\fU_w(K)\T\fU_{w'}(K)$ into $\fU_{w''}(K)$ where 
$w''$ is the product of $w,w'$ in $\fU(\{1\})=W$};

(d) {\it for $w_1,w_2,w'_1,w'_2$ in $W$, the multiplication in $\fG(K)$ 
maps $$\fG_{w_1,-w'_1}\T\fG_{w_2,-w'_2}(K)$$ into $\fG_{w_3,-w'_3}(K)$ where $(w_3,w'_3)$ is the product of 
$(w_1,-w'_1),(w_2,-w'_2)$ in \lb $\fG(\{1\})=W\T W$.}

\subhead 2.12\endsubhead
We now assume that $K$ is as in 1.1(i). The map $\fU_w(K)\T\fU_{w'}(K)@>>>\fU_{w''}(K)$ in 2.11(c) can be viewed
as a map $U_w^+(K)\T U_{w'}^+(K)@>>>U_{w''}^+(K)$ (multiplication in $U^+(K)$); using the definitions we see 
that this map
is admissible. Hence the map $\mu:U^+(K)\T U^+(K)@>>>U^+(K)$ (multiplication in $U^+(K)$) is admissible. The map 
$\fG_{w_1,-w'_1}\T\fG_{w_2,-w'_2}(K)@>>>\fG_{w_3,-w'_3}(K)$ in 2.11(d) can be viewed as a map 
$G_{w_1,-w'_1}\T G_{w_2,-w'_2}(K)@>>>G_{w_3,-w'_3}(K)$ (multiplication in $U(K)$); using the definitions we see 
that this map is admissible. Hence the map $\mu':G(K)\T G(K)@>>>G(K)$ (multiplication in $G(K)$) is admissible.
Applying these results with $K$ replaced by $K'=\RR(t)_{>0}$ we see that $\mu$ and $\mu'$ 
induce admissible maps of non-negative $\cz$-structures $\un\mu:\un{U^+(K')}\T\un{U^+(K')}@>>>\un{U^+(K')}$ and
$\un\mu':\un{G(K')}\T\un{G(K')}@>>>\un{G(K')}$, that is, $\un\mu:U^+(\cz)\T U^+(\cz)@>>>U^+(\cz)$ and
$\un\mu:G(\cz)\T G(\cz)@>>>G(\cz)$, which define monoid structures on $U^+(\cz)$ and $G(\cz)$. 

From the definition there is a unique homomorphism of monoids $c:\fU(\cz)@>>>U^+(\cz)$ which takes 
$i^a$ to the zone of $U^+(K')$ containing $x_i(t^a)$ (for $i\in I,a\in\cz$) and a unique homomorphism of monoids 
$c':\fG(\cz)@>>>G(\cz)$  which takes $i^a$ to the zone of $G_{s_i,-1}(K')$ 
containing $x_i(t^a)$, $(-i)^a$ to the zone of of $G_{1,-s_i}(K')$ containing $y_i(t^a)$ and
$\ui^a$ to the zone of $G_{1,-1}(K')$ containing $i(t^a)$ (for $i\in I,a\in\ZZ$).

\subhead 2.13\endsubhead
Assume that $K$ is as in 1.1(i)-(iii). Let $w\in W$ be such that $|w|=m$. For any $\ii\in\co_w$ we define 
$e_\ii:K^m@>>>\fU_w(K)$ by the same formula as in 2.9. We state the following property.

(a) {\it $e_\ii:K^m@>>>\fU_w(K)$ is a bijection.}
\nl
Let $w,w'$ in $W$ be such that $|w|+|w'|+r=M$. For any $\hh\in\co_{w,-w'}$ we define 
$\th_\hh:K^M@>>>\fG_{w,-w'}(K)$ by the same formula as in 2.10. We state the following property.

(b) {\it $\th_\hh:K^M@>>>\fG_{w,-w'}(K)$ is a bijection.}
\nl
Now (a),(b) are already known when $K$ is as in 1.1(i) and are obvious whe $K=\{1\}$. Moreover the maps
in (a),(b) are surjective. It remains to show that the maps in (a),(b) are injective when $K=\cz$.
It is then enough to show that the map in (a) (resp.(b)) for $K=\cz$ is injective after composition with 
$c$ (resp. $c'$) in 2.12. But that composition is a bijection since it is obtained by applying the functor of
taking zones to the corresponding map for $K=K'$ which is already known to be a bijection. This proves (a),(b).

This argument shows also that $c$ defines for any $w\in W$ a bijection $\fU_w(\cz)@>>>U^+_w(\cz)$ and that
$c'$ defines for any $w,w'$ in $W$ a bijection $\fG_{w,-w'}(\cz)@>>>G_{w,-w'}(\cz)$. It follows that

(c) {\it $c:\fU(\cz)@>>>U^+(\cz)$ is an isomorphism and $c':\fG(\cz)@>>>G(\cz)$ is an isomorphism.}
\nl
Thus, for any $K$ as in 1.1(i)-(iii) we have 

(d) $\fU(K)=U^+(K)$, $\fG(K)=G(K)$;
\nl
for $K$ as in 1.1(i) this is already known from 2.9(d) and 2.10(f). 

For any $K$ as in 1.1(i)-(iii) we have the following results.

(e) {\it Let $w\in W$. The bijections $e_\ii$ in (a) define a positive $K$-structure on $\fU_w(K)$. These 
positive $K$-structures can be viewed as a non-negative $K$-structure on $\fU(K)$ with pieces 
$\fU_w(K),w\in W$. The multiplication map $\fU(K)\T\fU(K)@>>>\fU(K)$ is admissible.}

(f) {\it Let $w,w'$ in $W$. The bijections $\th_\hh$ in (b) define a positive $K$-structure on $\fG_{w,-w'}(K)$.
These positive $K$-structures can be viewed as a non-negative $K$-structure on $\fG(K)$ with pieces
$\fG_{w,-w'}(K),(w,w')\in W\T W$. The multiplication map $\fG(K)\T\fG(K)@>>>\fG(K)$ is admissible.}
\nl
When $K$ is as in 1.1(i) these statements are already known, see 2.7, 2.8 and 2.9(d), 2.10(f). When $K=\{1\}$ 
these statements are trivial. When $K=\cz$ these statements can be deduced from the corresponding statements with
$K$ replaced by $\RR(t)_{>0}$ using (c).

\subhead 2.14\endsubhead
Let $\cn=\{0,1,2,\do\}\sub\cz$.
Let $w\in W$ with $|w|=m$. For $\ii\in\co_w$ let $\fU_w(\cn)$ be the image of the map $e_\ii:\cn^m@>>>\fU(\cz)$ 
given by $(a_1,a_2,\do,a_m)\m i_1^{a_1}i_2^{a_2}\do i_m^{a_m}$; this image is independent of the choice of $\ii$.
To prove the last statement it is enough to show that if $\ii,\ii'$ in $\co_w$ are connected by a 
single braid move then $e_\ii(\cn^m)=e_{\ii'}(\cn^m)$. This is a property of $R(i,j)$ in 2.9(ii). 
(The computation of $R(i,j)$ is reduced in 2.3 to the case where $s_is_j$ has order $2$ or $3$. When this order
is $3$ we use that $a+b-\min(a,c)\in\cn$ if $a,b,c\in\cn$.) Let $\fU(\cn)=\sqc_{w\in W}\fU_w(\cn)$. This is the 
submonoid of $\fU(\cz)$ generated by $i^a$ with $i\in I,a\in\cn$.

\subhead 2.15\endsubhead
Let $w,w'$ in $W$ be such that $|w|+|w'|+r=M$. For $\hh\in\co_{w,-w'}$ let 
$(\cn^M)'=(a_1,a_2,\do\a_m)\in\cn^M;a_k=0 \text{ whenever }h_k=\ui\text{ for some }i\}$ and let 
$\fG_{w,-w'}(\cn)$ be the image of the map $(\cn^M)'@>>>\fG(\cz)$ given by 
$(a_1,a_2,\do,a_m)\m h_1^{a_1}h_2^{a_2}\do h_M^{a_M}$; this image is independent of the choice of $\hh$. (We use 
the argument in 2.14 and the fact that, if $a,b,c$ are in $\cn$ and $b=0$ then $c-\min(a+c,b)=c\in\cn$, 
$\min(a+c,b)-b=0$.) Let $\fG(\cn)=\sqc_{w,w'\text{ in }W}\fG_{w,-w'}(\cn)$. This is the submonoid of $\fG(\cz)$ 
generated by $i^a$ and $(-i)^a$ with $i\in I,a\in\cn$. We have an isomorphism of monoids 
$\fU(\cn)\T\fU(\cn)@>\si>>\fG(\cn)$ given by $(i^a,1)\m i^a$, $(1,i^a)\m(-i)^a$.

\subhead 2.16\endsubhead
Let $K$ be as in 1.1(i)-(iii). Let $w\in W,m=|w|$. For $i\in I,a\in K$ such that $|s_iw|=m-1$ we define 
$T_{i,a}:\fU_w(K)@>>>\fU_w(K)$ as follows. Let $x\in\fU_{w_I}(K)$. We choose $\ii=(i_2,\do,i_m)\in\co_{s_iw}$ 
and we set $T_{i,a;\ii}(x)=i^{a_1a}i_2^{a_2}\do i_m^{a_m}$ where $x=i^{a_1}i_2^{a_2}\do i_m^{a_m}$ and 
$(a_1,a_2,\do,a_m)\in K^m$ is uniquely determined by $x$. We show that $T_{i,a;\ii}(x)$ is independent of the 
choice of $\ii$. Using Iwahori's lemma, we see that it is enough to show that $T_{i,a;\ii}(x)=T_{i,a;\ii'}(x)$ 
if $\ii,\ii'$ in $\co_{s_iw}$ are connected by a single braid move. In this case the desired result follows 
from 2.9(ii). We set $T_{i,a}(x)=T_{i,a;\ii}(x)$ where $\ii$ is any sequence in $\co_{s_iw}$. This defines the 
map $T_{i,a}$. We have $T_{i,a}T_{i,a'}=T_{i,aa'}$ for any and $a,a'$ in $K$.

\subhead 2.17\endsubhead
Let $K$ be as in 1.1(i)-(iii). Let $\l\in\cx,w\in W,m=|w|$. We will define a map
$$\th_\l:\fU_w(K)@>>>K.\tag a$$
For any $\ii=(i_1,i_2,\do,i_m)\in\co_w$ we define $\th_{\l,\ii}:\fU_w(K)@>>>K$ by
$$i_1^{a_1}i_2^{a_2}\do i_m^{a_m}\m a_1^{c_1}a_2^{c_2}\do a_m^{c_m}$$ 
where $(a_1,a_2,\do,a_m)\in K^\nu$ and $c_k=\la s_{i_1}s_{i_2}\do s_{i_{k-1}}(i_k),\l\ra\in\ZZ$ for 
$k\in[1,m]$. We show that $\th_{\l,\ii}=\th_{\l,\ii'}$ for any $\ii,\ii'$ in $\co_w$. It is enough to show 
this assuming that $\ii,\ii'$ are related by a single braid move.

Assume first that for some $l$ we have $i_l=i,i_{l+1}=j$ where $\la i,j^*\ra=0$ and
$a_l=a,a_{l+1}=b$. We set $y=s_{i_1}s_{i_2}\do s_{i_{l-1}}\in W$. It is enough to show that
$$a^{\la y(i),\l\ra}b^{\la ys_i(j),\l\ra}=b^{\la y(j),\l\ra}a^{\la ys_j(i),\l\ra}.$$
This follows from $ s_i(j)=j,s_j(i)=i$.

Assume next that for some $l$ we have $i_l=i,i_{l+1}=j,i_{l+2}=i$ where $\la i,j^*\ra=\la j,i^*\ra=-1$ and
$a_l=a,a_{l+1}=b,a_{l+2}=c$. We set $y=s_{i_1}s_{i_2}\do s_{i_{l-1}}\in W$. It is enough to show that
$$\align&a^{\la y(i),\l\ra}b^{\la ys_i(j),\l\ra}c^{\la ys_is_j(i),\l\ra}\\&=
(bc/(a+c))^{\la y(j),\l\ra}(a+c)^{\la ys_j(i),\l\ra}(ab/(a+c)^{\la ys_js_i(j),\l\ra}.\endalign$$
This follows from $s_js_i(j)=i$, $s_is_j(i)=j$, $s_i(j)=i+j=s_j(i)$.

Assume next that for some $l$ we have $i_l=i,i_{l+1}=j,i_{l+2}=i,i_{l+3}=j$ where $\la i,j^*\ra=-2$ and
$a_l=a,a_{l+1}=b,a_{l+2}=c,a_{l+3}=d$. We set $y=s_{i_1}s_{i_2}\do s_{i_{l-1}}\in W$. It is enough to show that
$$\align&a^{\la y(i),\l\ra}b^{\la ys_i(j),\l\ra}c^{\la ys_is_j(i),\l\ra}d^{\la ys_is_js_i(j),\l\ra}\\&=
(bc^2d/E)^{\la y(j),\l\ra}(E/A)^{\la ys_j(i),\l\ra}(A^2/E)^{\la ys_js_i(j),\l\ra}
(abc/A)^{\la ys_js_is_j(i),\l\ra}\endalign$$
where $A,E$ are as in 2.4(d). This follows from
$s_i(j)=i+j$, $s_j(i)=i+2j$, $s_is_j(i)=i+2j$, $s_js_i(j)=i+j$, $s_is_js_i(j)=j$, $s_js_is_j(i)=s_i$.

Finally we assume that for some $l$ we have $i_l=i,i_{l+1}=j,i_{l+2}=i,i_{l+3}=j,i_{l+4}=i,i_{l+5}=j$
where $\la i,j^*\ra=-3$. In this case an argument similar to the one above applies.
We see that  $\th_{\l,\ii}$ is indeed independent of the choice of $\ii$; we denote it by $\th_\l$.
Note that for $\l,\l'$ in $\cx$ we have $\th_{\l+\l'}(u)=\th_\l(u)\th_{\l'}(u)$ for any $u\in\fU_w(K)$.

In the case where $K=\cz,w=w_I$, a map like $\th_\l$ (restricted to $\fU_{w_I}(\cn)$) is defined in 
\cite{L90, 2.9(a)}.

In the case where $\l\in\cx^+$ and $K,\kk$ are as in 1.1, we define a map $\ti\th_\l:G@>>>\kk$, $g\m x$
where $g\in G$ and $x\in\kk$ is the coefficient of the highest weight vector in the canonical basis of $\L_\l$
in the vector obtained by applying $g$ to the lowest weight weight vector in the canonical basis of $\L_\l$.
One can show that in this case (with $w=w_I$):

(b) {\it $\th_\l$ is the composition $\fU_{w_I}(K)@>\si>>U^+_{w_I}(K)\sub G@>\ti\th_\l>>\kk$.}

\head 3. The positive part of a flag manifold\endhead
\subhead 3.1\endsubhead
We assume that $K,\kk$ are as in 1.1(i).
 Let $M=\nu+r$. Let 
$$\co'=\sqc_{w\in W}\co_{w,-w_Iw},$$
a set of sequences of length 
$M$ in $I\sqc(-I)\sqc\uI$, see 2.8. Let $\G'$ be the set of all pairs $(\hh,\aa)\in\co'\T K^M$. With notation 
of 2.8 we have $\G'=\sqc_{w\in W}\G_{w,-w_Iw}$. We define a map $\r:\G'@>>>G/U^-$ by the requirement that the 
restriction of $\r$ to $\G_{w,-w_Iw}$ is 
$$((h_1,h_2,\do,h_M),(a_1,a_2,\do,a_M))\m h_1^{a_1}h_2^{a_2}\do h_M^{a_M}(w\i w_I)\dot{}U^-$$
where $h_s^{a_s}=x_i(a_s)$ if $h_s=i$, $h_s^{a_s}=y_i(a_s)$ if $h_s=-i$ and $h_s^{a_s}=i(a_s)$ if $h_s=\ui$.

We view $\G'$ as the set of vertices of a graph with two kind of edges: internal edges and external edges. An 
internal edge is an edge joining two elements in the same $\G_{w,-w_Iw}$ for the graph structure of 
$\G_{w,-w_Iw}$ described in 2.8. An external edge is one joining $((h_1,h_2,\do,h_M),(a_1,a_2,\do,a_M))\in\G'$, 
$((h'_1,h'_2,\do,h'_M),(a'_1,a'_2,\do,a'_M))\in\G'$ where

(a) $(h_{M-1},h_M)=(i,\ui)$, $(h'_{M-1},h'_M)=(-i,\ui)$ for some $i\in I$, 

$(a'_{M-1},a'_M)=(a_{M-1}\i,a_{M-1}a_M\i)$,

$h'_s=h_s,a'_s=a_s$ for $s\n\{M-1,M\}$;
\nl
this edge joins an element of $\G_{ys_i,-w_Iys_i}$ with an element of $\G_{y,-w_Iy}$ for some $y\in W$ such that 
$|ys_i|>|y|$.

We regard $\co'$ as the set of vertices of a graph in which $\hh,\hh'$ are joined if either $\hh,\hh'$ belong 
to the same $\co_{w,-w_Iw}$ and are joined in the graph $\co_{w,-w_Iw}$ (see 2.8), or 

(b) $\hh,\hh'$ are of the form  $\hh=(h_1,\do,h_{M-2},i,\ui)$, $\hh'=(h_1,\do,h_{M-2},-i,\ui)$ with $i\in I$. 
\nl
The map $\G'@>>>\co'$, $(\hh,\aa)\m\hh$, respects the graph structures. 

\subhead 3.2\endsubhead
For any $\hh\in\co'$ we define $\r_\hh:K^M@>>>G/U^-$ by $\r_\hh(\aa)=\r(\hh,\aa)$. We show:

(a) {\it The image of $\r_\hh$ is independent of the choice of $\hh\in\co'$; we denote this image by 
$(G/U^-)(K)_{>0}$. For any $\hh\in\co'$, $\r_\hh$ defines a bijection $$\r'_\hh:K^M@>\si>>(G/U^-)(K)_{>0}.$$ The 
bijections $\r'_\hh (\hh\in\co')$ define a positive $K$-structure on $(G/U^-)(K)_{>0}$.}
\nl
Assume first that $\hh,\hh'$ in $\co'$ are joined in the graph $\co'$. Then 

(b) {\it $\r_\hh,\r_{\hh'}$ are related by $\r_\hh=\r_{\hh'}\s$ where $\s:K^M@>>>K^M$ is a bi-admissible 
bijection. In particular, the image of $\r_\hh$ is equal to the image of $\r_{\hh'}$ and $\r_\hh$ is injective 
if and only if $\r_{\hh'}$ is injective.}
\nl
If $\hh,\hh'$ are in the same $\co_{w,-w_Iw}$, this follows from results in 2.8. Assume now that 
$\hh,\hh'$ are as in 3.1(b) so that $\hh\in\co_{ys_i,-w_Iys_i}$, $\hh'\in\co_{y,-w_Iy}$ for some $y\in W$ 
such that $|ys_i|>|y|$. Then (b) follows from the calculation
$$\align&x_i(a)i(b)(s_i\i y\i w_I)\dot{}U^-=y_i(a\i)i(ab\i)\ds_iy_i(a\i b^2)(s_i\i y\i w_I)\dot{}U^-\\&=
y_i(a\i)i(ab\i)\ds_i(s_i\i y\i w_I)\dot{}U^-=y_i(a\i)i(ab\i)(y\i w_I)\dot{}U^-\endalign$$
for $a,b$ in $K$.
(We have used 2.1(a) and $y_i(c)\dw\in\dw U^-$, $\ds_i\dw=(s_iw)\dot{}$ for $c\in\kk$ and 
$w\in W$ with $|s_iw|>|w|$, or more precisely, $w=s_i\i y\i w_I$.) In this case $\s:K^M@>>>K^M$ is 
$(a_1,\do,a_{M-2},a,b)\m(a_1,\do,a_{M-2},a\i,ab\i)$. In particular, $\s$ is a bi-admissible  bijection.

\mpb

Using (b), we see that to prove (a) it is enough to show (c),(d) below.

(c) {\it For any $\hh\in\co_{w_I,-1}$, $\r_\hh:K^M@>>>G/U^-$ is injective.}

(d) {\it the graph $\co'$ is connected.}
\nl
We prove (c). Let $\hh$ in $\co_{w_I,-1}$. Let $\aa\in K^M$, $\aa'\in K^M$ be such 
that $\r_\hh(\aa)=\r_\hh(\aa')$, that is $$\g_{w_I,-1}(\aa)U^-=\g_{w_I,-1}(\aa')U^-$$ (notation of 2.8). Since
$\g_{w_I,-1}(\aa)\in B^+$, $\g_{w_I,-1}(\aa')\in B^+$ and $B^+\cap U^-=\{1\}$, it follows that
$$\g_{w_I,-1}(\aa)=\g_{w_I,-1}(\aa').$$ Using the injectivity of $\g_{w_I,-1}:K^M@>>>G$
(see 2.8) we deduce that $\aa=\aa'$. This proves (c).

We prove (d). As mentioned in 2.8, the graph $\co_{w,-w_Iw}$ is connected for any $w\in W$.
In particular $\co_{w_I,-1}$ is contained in a connected component $\co'_0$ of $\co'$.
We show by descending induction on $|w|$ that $\co_{w,-w_Iw}$ is contained in $\co'_0$.
We can assume that $|w|<\nu$. We can find $i\in I$ such that $|ws_i|>|w|$. 
Since the graph  $\co_{w,-w_Iw}$ is connected, we can find $\hh'\in\co_{w,-w_Iw}$
such that $\hh'=(h_1,h_2,\do,h_{M-2},-i,\ui)$. 
Then $\hh=(h_1,h_2,\do,h_{M-2},i,\ui)$ is joined with $\hh'$ in $\co'$ and is in $\co_{ws_i,w_Iws_i}$.
By the induction hypothesis we have $\hh\in\co'_0$ hence $\hh'\in\co'_0$ so that 
$\co_{w,-w_Iw}\sub\co'_0$. This completes the induction. We see that $\co'=\co'_0$ and (d) holds.
This completes the proof of (a).

\mpb

We show:

(e) {\it If $g_1\in G(K)$ and $gU^-\in(G/U^-)(K)_{>0}$ then $g_1gU^-\in(G/U^-)(K)_{>0}$. Hence there is a well 
defined action of the monoid $G(K)$ on $(G/U^-)(K)_{>0}$ given by $g_1:gU^-\m g_1gU^-$.}
\nl
Assume first that $g_1=x_i(a)$ with $i\in I,a\in K$.
We can find $\hh\in\co_{w_I,-1}$ such that $\hh$ starts with $i$ and $\aa=(a_1,a_2,\do,a_M)\in K^M$ such that
$gU^-=\r_\hh(\aa)$. We have $x_i(a)gU^-=x_i(a)\r_\hh(\aa)=\r_\hh(a_1+a,a_2,a_3,\do,a_M)\in(G/U^-)(K)_{>0}$, as
required. Next we assume that $g_1=y_i(a)$ with $i\in I,a\in K$.
We can find $\hh\in\co_{1,-w_I}$ such that $\hh$ starts with $-i$ and $\aa=(a_1,a_2,\do,a_M)\in K^M$ such that
$gU^-=\r_\hh(\aa)$. We have $y_i(a)gU^-=y_i(a)\r_\hh(\aa)=\r_\hh(a_1+a,a_2,a_3,\do,a_M)\in(G/U^-)(K)_{>0}$, as
required. We now assume that $g_1=i(a)$ with $i\in I,a\in K$.
We can find $\hh\in\co_{w_I,-1}$ such that $\hh$ starts with $\ui$ and $\aa=(a_1,a_2,\do,a_M)\in K^M$ such that
$gU^-=\r_\hh(\aa)$. We have $i(a)gU^-=i(a)\r_\hh(\aa)=\r_\hh(aa_1,a_2,a_3,\do,a_M)\in(G/U^-)(K)_{>0}$, as 
required.
It remains to use that the monoid $G(K)$ is generated by $x_i(a),y_i(a),i(a)$ for various $i\in I,a\in K$.

\mpb

From the definitions we see that for $w,w'$ in $W$, the map 
$$G_{w,-w'}(K)\T (G/U^-)(K)_{>0}@>>>(G/U^-)(K)_{>0}$$ (restriction of the $G(K)$-action
in (e)) is admissible. Hence the map 

(f) $G(K)\T(G/U^-)(K)_{>0}@>>>(G/U^-)(K)_{>0}$ 
\nl
given by the $G(K)$-action is admissible. 

Now let $K'=\RR(t)_{>0}$. We set $(G/U^-)(\cz)_{>0}=\un{(G/U^-)(K')_{>0}}$. By passage to zones in (f) we obtain
an admissible map $G(\cz)\T(G/U^-)(\cz)_{>0}@>>>(G/U^-)(\cz)_{>0}$ which is an action of $G(\cz)$.

\subhead 3.3\endsubhead
Now $T$ acts on $G/U^-$ by $t:gU^-\m gtU^-$. We show:

(a) {\it This restricts to an action of the group $T(K)$ on $(G/U^-)(K)_{>0}$.}
\nl
Let $gU^-\in(G/U^-)(K)_{>0}$, $t\in T(K)$. Let $(i_1,i_2,\do,i_\nu)\in\co_{w_I}$ and let
$j_1,\do,j_r$ be a list of the elements in $I$. We have $t=j_1(b_1)\do j_r(b_r)$ with $(b_1,\do,b_r)\in K^r$.
We can find $(a_1,a_2,\do,a_{\nu+r})\in K^{\nu+r}$
such that 
$$gU^-=x_{i_1}(a_1)\do x_{i_\nu}(a_\nu)j_1(a_{\nu+1})\do j_r(a_{\nu+r})U^-.$$
We have 
$$gtU^-=x_{i_1}(a_1)\do x_{i_\nu}(a_\nu)j_1(a_{\nu+1}b_1)\do j_r(a_{\nu+r}b_r)U^-.$$
Hence $gtU^-\in(G/U^-)(K)_{>0}$.

We show: 

(b) {\it The obvious map $(G/U^-)(K)_{>0}/T(K)@>>>(G/U^-)/T$ is injective.}
\nl
Let $gU^-\in(G/U^-)(K)_{>0}$, $g'U^-\in(G/U^-)(K)_{>0}$, $t\in T$ be such that $gtU^-=g'U^-$. We must show that
$t\in T(K)$. With notation in the proof of (a) we have

$gU^-=x_{i_1}(a_1)\do x_{i_\nu}(a_\nu)j_1(a_{\nu+1})\do j_r(\a_{\nu+r})U^-$,

$g'U^-=x_{i_1}(a'_1)\do x_{i_\nu}(a'_\nu)j_1(a'_{\nu+1})\do j_r(a'_{\nu+r})U^-$,

$t=j_1(c_1)\do j_r(c_r)$,
\nl
with $(c_1,\do,c_r)\in(\kk^*)^r$, 
$(a_1,a_2,\do,a_{\nu+r})\in K^{\nu+r}$, $(a'_1,a'_2,\do,a'_{\nu+r})\in K^{\nu+r}$.
\nl
Our assumption implies
$$\align&x_{i_1}(a'_1)\do x_{i_\nu}(a'_\nu)j_1(a'_{\nu+1})\do j_r(a'_{\nu+r})U^-\\&=
x_{i_1}(a_1)\do x_{i_\nu}(a_\nu)j_1(a_{\nu+1}c_1)\do j_r(a_{\nu+r}c_r)U^-.\endalign$$
It follows that $B^+j_1(a_{\nu+1}c_1)\do j_r(a_{\nu+r}c_r)U^-=B^+j_1(a'_{\nu+1})\do j_r(a'_{\nu+r})U^-$.
By a property of Bruhat decomposition we deduce that
$$j_1(a_{\nu+1}c_1)\do j_r(a_{\nu+r}c_r)=j_1(a'_{\nu+1})\do j_r(a'_{\nu+r})$$
hence $a_{\nu+1}c_1=a'_{\nu+1},\do,a_{\nu+r}c_r=a'_{\nu+r}$, so that $(c_1,\do,c_r)\in K^r$ and $t\in T(K)$. This 
proves (b).

\subhead 3.4\endsubhead
We assume that $K,\kk$ are as in 1.1(i).
Let $\cb$ the set of subgroups of $G$ of the form $gB^-g\i$ for some $g\in G$. Following \cite{L94} we set
$$\align&\cb(K)_{>0}=\{B\in\cb;B=uB^+u\i\text{ for some } u\in U^-_{w_I}(K)\}\\&
=\{B\in\cb;B=u'B^-u'{}\i\text{ for some } u'\in U^+_{w_I}(K)\}.\endalign$$
(The last equality is proved in \cite{L94, 8.7}.) 

Consider the set of orbits for the free action of $T$ on $G/U^-$ given by $t:gU^-\m gtU^-$. 
This set can be identified with $\cb$ by $gU^-\m gB^-g\i$. Under this identification, 
$\cb(K)_{>0}$ becomes the set of orbits for the free action of $T(K)$ on $(G/U^-)(K)_{>0}$ (the
restriction of $T$-action above); here we use 3.3(b).

For any $\hh=(h_1,h_2,\do,h_M)\in\co'$, 
there is a unique free action $t:\hat\aa\m t\cir_\hh\hat\aa$ of $T(K)$ on $K^M$ 
such that the bijection $K^M@>>>(G/U^-)(K)_{>0}$ given by $(a_1,a_2,\do,a_M)\m\r_\hh(a_1,a_2,\do,a_M)$ is 
compatible with the $T(K)$-actions on the two sides. Let $k_1<k_2<\do<k_\nu$ be the subsequence of 
$1,2,\do,M$ such that $h_{k_1},h_{k_2},\do,h_{k_\nu}$ are in $I\sqc(-I)$. We define an 
imbedding $\e'_\hh:K^\nu@>>>K^M$ by $(a_1,a_2,\do,a_\nu)\m(a'_1,a'_2,\do,a'_M)$ where
$a'_{k_j}=a_j$ for $j=1,\do,\nu$ and $a'_e=1$ if $e\in\{1,2,\do,M\}-\{k_1,k_2,\do,k_\nu\}$.
To any $T(K)$-orbit in $K^M$ for the $\cir_\hh$-action we associate the unique element $\aa\in K^\nu$ 
such that our $T(K)$-orbit contains $\e'_\hh(\aa)$. This defines a map $\bar\d'_\hh:K^M/T(K)@>>>K^\nu$ 
which is easily seen to be a bijection. Let $\d'_\hh:K^M@>>>K^\nu$ be the composition
of $\bar\d'_\hh$ with the obvious map $K^M@>>>K^M/T(K)$. There is a unique map 
 $\bar\t'_\hh:K^\nu@>>>\cb(K)_{>0}$ such that the following diagram is commutative:
$$\CD
K^M@>\t'_\hh>>(G/U^-)(K)_{>0}\\$$
@V\d'_\hh VV @VVV\\
K^\nu@>\bar\t'_\hh>>\cb(K)_{>0}
\endCD$$
Here $\t'_\hh$ is the bijection given by $\aa\m\r_\hh(\aa)$ and the right vertical map 
is the obvious orbit map. Clearly, $\bar\t'_\hh$ is a bijection. We show:

(a) {\it The bijections $K^\nu@>\bar\t'_\hh>>\cb(K)_{>0}$ for various $\hh\in\co'$ define a positive $K$-structure
on $\cb(K)_{>0}$.}
\nl
Let $\hh,\hh'$ be elements of $\co'$. Let 
$$A=\t'_{\hh'}{}\i\t'_\hh:K^M@>>>K^M,\qua \bar A=\bar\t'_{\hh'}{}\i\bar\t'_\hh:K^\nu@>>>K^\nu.$$
 We have a commutative diagram
$$\CD
K^M@>A>>K^M\\$$
@V\d'_\hh VV  @V\d'_{\hh'}VV\\
K^\nu@>\bar A>>K^\nu
\endCD$$
From 3.2(a) it follows that $A$ is admissible. From the definitions we see that the 
surjective maps
$\d'_\hh,\d'_{\hh'}$ are admissible. Moreover $\d'_\hh\e'_\hh=1$. Thus we have $\bar A=\d'_{\hh'}A\e'_\hh$. It 
follows that $\bar A$ is admissible. 
Interchanging  the roles of $\hh,\hh'$ we see that $\bar A\i$ is admissible. This proves (a).

\mpb

From 3.2(a) and the definitions we deduce:

(b) {\it There is a well defined action of the monoid $G(K)$ on $\cb(K)_{>0}$ given by 
$g_1:gB^-g\i\m g_1gB^-g\i g_1\i$ where $gB^-g\i\in \cb(K)_{>0}$.}

\mpb

From the definitions we see that for $w,w'$ in $W$, the map 
$$G_{w,-w'}(K)\T\cb(K)_{>0}@>>>\cb(K)_{>0}$$
 (restriction of the $G(K)$-action in (b)) is admissible. Hence the map 

(c) $G(K)\T\cb(K)_{>0}@>>>\cb(K)_{>0}$ 
\nl
given by the $G(K)$-action is admissible. 

Now let $K'=\RR(t)_{>0}$. We set $\cb(\cz)_{>0}=\un{\cb(K')_{>0}}$. By passage to zones in (c) we obtain

(d) an admissible map $G(\cz)\T\cb(\cz)_{>0}@>>>\cb(\cz)_{>0}$ which is an action of $G(\cz)$.

\head 4. The involutions $\Ph,\Ph'$ of $\fU_{w_I}(K)$\endhead
\subhead 4.1\endsubhead
We assume that $K,\kk$ are as in 1.1(i).
Let $\ii=(i_1,i_2,\do,i_\nu)\in\co_{w_I}$. We define $\x_\ii:K^\nu@>\si>>\cb(K)_{>0}$ by 
$(a_1,a_2,\do,a_\nu)\m gB^-g\i$ where 
$$g=x_{i_1}(a_1)x_{i_2}(a_2)\do x_{i_\nu}(a_\nu)\in G.$$ 
We define $\x'_\ii:K^\nu@>\si>>\cb(K)_{>0}$ by $(a_1,a_2,\do,a_\nu)\m g'B^+g'{}\i$ where 
$$g'=y_{i_1}(a_1)y_{i_2}(a_2)\do y_{i_\nu}(a_\nu)\in G.$$ 
Clearly,

(a) {\it the bijections $\x_\ii,\x'_\ii$ are parts of the positive $K$-structure of $\cb(K)_{>0}$ described 
in 3.4.}

\mpb

For $u\in U^+_{w_I}(K)$ 
we have by 2.5(d) $u=\ph'(u)\dw_Ib$ where $b\in B^-$ and $\ph'$ is a well defined map
$U^+_{w_I}(K)@>>>U^-$. We show:

(b) $u\m\ph'(u)$ is a bijection $U^+_{w_I}(K)@>>>U^-_{w_I}(K)$.
\nl
The map $U^+_{w_I}(K)@>>>\cb(K)_{>0}$, $u\m uB^-u\i$, is a bijection; the map $U^-_{w_I}(K)@>>>\cb(K)_{>0}$
$u\m uB^+u\i$, is a bijection. Hence if $u\in U^+_{w_I}(K)$, there is a unique $u'\in U^-_{w_I}(K)$ such that $uB^-u\i=u'B^+u'{}\i$ so 
that 
$$uB^-u\i=\ph'(u)\dw_IB^-\dw_I\i\ph'(u)\i=\ph'(u)B^+\ph'(u)\i=u'B^+u'{}\i$$ 
and $u'=\ph'(u)$. Moreover 
$u\m u'$ is a bijection $U^+_{w_I}(K)@>>>U^-_{w_I}(K)$ so that $u\m\ph'(u)$ is a bijection 
$\ph':U^+_{w_I}(K)@>>>U^-_{w_I}(K)$. This proves (b).

\mpb

For $u\in U^+_{w_I}(K)$ we have by 2.5(c) $u\in U^-\dw_IB^-$ hence $u\i\in B^-\dw_IU^-$. Thus we have 
$u\i=\ph(u)\i\dw_Ib$ where $b\in B^-$ and $\ph$ is a well defined map $U^+_{w_I}(K)@>>>U^-$. We show:

(c)  $u\m\ph(u)$ is a bijection $U^+_{w_I}(K)@>>>U^-_{w_I}(K)$.
\nl
From the definition of $\ph'$, for $u\in U^+_{w_I}(K)$ we have
$\ph'{}\i(\Ps(u))=\Ps(u)\dw_Ib$ for some $b\in B^-$ (with $\Ps$ as in 2.1.) Hence
$\Ps(u)=\ph'{}\i(\Ps(u))\dw_I b_1$ for some $b_1\in B^+$ and
$u=b_2\dw_I\Ps(\ph'{}\i(\Ps(u)))$ for some $b_2\in B^-$.

From the definition we have $\ph(u)\i=u\i\dw_Ib_3$, for some $b_3\in B^+$ hence
$u=\dw_Ib_3\ph(u)=b_4\dw_I\ph(u)$ for some $b_4\in B^-$. Thus 
$b_4\dw_I\ph(u)=b_2\dw_I\Ps(\ph'{}\i(\Ps(u)))$. From properties of Bruhat decomposition we deduce
$\ph(u)=\Ps(\ph'{}\i(\Ps(u)))$. Thus (c) follows from (b).

\mpb

From the definition, for $u\in U^+_{w_I}(K)$ we have $\ph(u)=b\dw_Iu$ for a uniquely defined $b\in B^+$.
Let $\ii=(i_1,i_2,\do,i_\nu)\in\co_{w_I}$. There is a unique bijection $z:K^\nu@>>>K^\nu$,
$(a_1,\do,a_\nu)\m(a'_1,\do,a'_\nu)$ such that 
$$\ph(x_{i_1}(a_1)\do x_{i_\nu}(a_\nu))=y_{i_1}(a'_1)\do y_{i_\nu}(a'_\nu)$$ 
where
$$y_{i_1}(a'_1)\do y_{i_\nu}(a'_\nu)=b\dw_Ix_{i_1}(a_1)\do x_{i_\nu}(a_\nu)$$
for some $b\in B^+$; moreover, by (a), $z$ is bi-admissible. 

Applying the involution $\io$ (see 2.1) we deduce
$$x_{i_1}(a'_1)\do x_{i_\nu}(a'_\nu)=b'\dw_Iy_{i_1}(a_1)\do y_{i_\nu}(a_\nu)$$
for some $b'\in B^-$ hence
$$y_{i_1}(a_1)\do y_{i_\nu}(a_\nu)=b''\dw_Ix_{i_1}(a'_1)\do x_{i_\nu}(a'_\nu)$$
for some $b''\in B^+$. It follows that $z^2=1$. 

Let 

(d) $e:\fU_{w_I}(K)@>\si>>U^+_{w_I}(K)$, $e':\fU_{w_I}(K)@>\si>>U^-_{w_I}(K)$
\nl
be the bijections induced by 2.9(d),(e). We define a bijection $\Ph:\fU_{w_I}(K)@>\si>>\fU_{w_I}(K)$ 
by $\Ph=e'{}\i\ph e$. We have $\Ph^2=1$ (we use the equality $z^2=1$).

We define a bijection $\Ph':\fU_{w_I}(K)@>\si>>\fU_{w_I}(K)$ by $\Ph'=e'{}\i\ph'e$.
Recall that the bijection $\ti\Ps:\fU_{w_I}(K)@>\si>>\fU_{w_I}(K)$ (see 2.9(f)) satisfies $\ti\Ph=e'{}\i\Ps e$
and $\ti\Ps^2=1$. We have 
$\Ph=\ti\Ps\Ph'{}\i\ti\Ps$ hence $1=\Ph^2=\ti\Ps\Ph'{}^{-2}\ti\Ps$ and $\Ps'{}^2=1$. It follows that

(e) $\Ph=\ti\Ps\Ph'\ti\Ps$.

\mpb

Now let $K'=\RR(t)_{>0}$. The bijection $\Ph:\fU_{w_I}(K')@>\si>>\fU_{w_I}(K')$ is bi-admissible since 
$z$ is so. Hence $\Ph'$ is also bi-admissible. Thus $\Ph$ induces by passage to zones a bi-admissible 
bijection $\fU_{w_I}(\cz)@>\si>>\fU_{w_I}(\cz)$ denoted  again by $\Ph$ (with $\Ph^2=1$), see also 
\cite{L97, 2.9}. Similarly, $\Ph'$ induces by passage to zones a bi-admissible bijection 
$\fU_{w_I}(\cz)@>\si>>\fU_{w_I}(\cz)$ denoted  again by $\Ph'$ (with $\Ph'{}^2=1$). Then (e) continues to
hold.

In type $A_2$ with $I=\{1,2\}$, $\Ph$ is given by
$$2^{a_2}1^{a_1}2^{a'_2}\m 2^{a_2/((a_2+a'_2)a'_2)}1^{(a_2+a'_2)/(a_1a_2)}2^{1/(a_2+a'_2)}.$$
In type $A_3$ with $I=\{1,2,3\}$ and $s_1s_3=s_3s_1$, $\Ph$ is given by
$$\align&2^{a_2}1^{a_1}3^{a_3}2^{a'_2}1^{a'_1}3^{a'_3}\m\\&
2^{\fra{a'_1a'_3}{a_1a_2a_3}}1^{\fra{a_1}{a'_1(a_1+a'_1)}}3^{\fra{a_3}{a'_3(a_3+a'_3)}}2^{\fra{(a_1+a'_1)(a_3+a'_3)}{a_1a_3a'_2}}
1^{\fra{1}{a_1+a'_1}}3^{\fra{1}{a_3+a'_3}}.\endalign$$

\subhead 4.2\endsubhead
Let $K$ be as in 1.1(i)-(iii). For $i\in I,a\in K$ we have

(a) $\Ph T_{i,a}=T_{i^!,a\i}\Ph:\fU_{w_I}(K)@>>>\fU_{w_I}(K)$
\nl
where $\Ph:\fU_{w_I}(K)@>>>\fU_{w_I}(K)$ is as in 4.1 (for $K=\{1\}$, $\Ph$ is the identity map) 
and $T_{i,a}$ is as in 2.16.
For $K$ as in 1.1(i) this follows from \cite{L97, 3.6}; for $K=\cz$ this follows from the case $K=\RR(t)_{>0}$ by passage to zones.
For $K=\{1\}$, (a) is trivial.

\mpb

For $K$ as in 1.1(i)-(iii) we conjecture the following result.

(b) $\th_\l\Ph=\th_{-\l^!}:\fU_{w_I}(K)@>>>K$ for any $\l\in\cx$ (notation of 2.17(a)); here $\l^!\in\cx$ is 
defined by $\la i,\l^!\ra=\la i^!,\l\ra$ for all $i\in I$.
\nl
When $K=\cz$, (b) can be deduced from \cite{L97, 4.9}. For general $K$ but for type $A_2$ with $I=\{i,j\}$, (b) 
follows from the identity:
$$a^{\la j,\l\ra}b^{\la i,\l\ra+\la j,\l\ra}c^{\la i,\l\ra}=
(\fra{a}{(a+c)c})^{-\la i,\l\ra}(\fra{a+c}{ab})^{-\la i,\l\ra-\la j,\l\ra}(\fra{1}{a+c})^{-\la j,\l\ra}$$
where $a,b,c$ are in $K$.

\mpb

When $K=\cz$ we define $h_0\in\fU_{w_I}(K)$ by $h_0=i_1^0i_2^0\do i_\nu^0$ where $\ii=(i_1,i_2,\do,i_\nu)$
is any element of $\co_{w_0}$. The following result appears in \cite{L97, 2.9}.

(c) {\it Assume that $K=\cz$. Then $\Ph:\fU_{w_I}(K)@>>>\fU_{w_I}(K)$ is the unique bijection satisfying (a)
for any $i\in I,a\in K$ and such that $\Ph(h_0)=h_0$.}

\subhead 4.3\endsubhead
Let $\l:I@>>>K$ be a map. We have the following results.

(a) {\it There is a unique monoid automorphism $S_\l:\fU(K)@>>>\fU(K)$ such that for any $i\in I,a\in K$ we have 
$S_\l(i^a)=i^{a\l(i)}$. This maps the subset $\fU_w(K)$ onto itself for any $w$ in $W$.}

(b) {\it There is a unique monoid automorphism $S_\l:\fG(K)@>>>\fG(K)$ such that for any $i\in I,a\in K$ we have 
$S_\l(i^a)=i^{a\l(i)}$, $S_\l((-i)^a)=(-i)^{a\l(i)\i}$, $S_f(\ui^a)=\ui^a$. This maps the subset $\fG_{w,-w'}(K)$ 
onto itself for any $w,w'$ in $W$.}
\nl
Indeed, one can check that the relations defining $\fU(K)$, $\fG(K)$ are respected by $S_\l$. Alternatively, 
assuming that $K,\kk$ are as in 1.1(i), let $\tS_\l:G@>>>G$ be given by conjugation by an element $t_\l\in T(K)$
 such that $i^*(t_\l)=\l(i)$ for any $i\in I$. It is enough to note that $\tS_\l$ maps $x_i(a)$ to $x_i(a\l(i))$ 
and $y_i(a)$ to $y_i(a\l(i)\i)$ for any $i\in I,a\in\kk$. The case where $K=\cz$ can be deduced from the case where 
$K=\RR(t)_{>0}$ by passage to zones. The case where $K=\{1\}$ is trivial.

\mpb

For $u\in U^+_{w_I}(K)$ with $K,\kk$ as in 1.1(i), we show:

(c) $\tS_\l(\ph(u))=\ph(\tS_\l(u))$.
\nl

We have $\ph(u)\in B^+\dw_Iu$, $\ph(t_\l ut_\l\i)\in B^+\dw_It_\l ut_\l\i$, hence 
$$t_\l\ph(u)t_\l\i\in B^+\dw_It_\l ut_\l\i.$$
 It follows that
$$\tS_\l(\ph(u))\in B^+\dw_I\tS_\l(u), \ph(\tS_\l(u))\in B^+\dw_I\tS_\l(u)$$ 
so that $\tS_\l(\ph(u))\ph(\tS_\l(u))\i\in B^+$. Since
$$\tS_\l(\ph(u))\ph(\tS_\l(u))\i\in U^-$$ and 
$B^+\cap U^-=\{1\}$, we deduce that 
$\tS_\l(\ph(u))\ph(\tS_\l(u))\i=1$ so that (c) holds.

\mpb

For $u\in\fU_{w_I}(K)$ we show: 

(d) $S_\l(\Ph(u))=\Ph(S_{\l\i}(u))$ where $\l\i:I@>>>K$ is given by $\l\i(i)=\l(i)\i$.
\nl
An equivalent statement is that $S_\l(e'{}\i(\ph(e(u))))=e'{}\i(\ph(e((S_{\l\i}(u))))$ where $e,e'$ are as in 
4.1(d). Let $u'=\ph(e(u))\in U^-_{w_I}(K)$. From the definitions we have 

$S_\l(e'{}\i(u'))=e'{}\i(\tS_{\l\i}(u'))$, $e((S_{\l\i}(u)))=\tS_{\l\i}(e(u))$. Hence it is enough to show 
that

$e'{}\i(\ti S_{\l\i}(u'))=e'{}\i(\ph(\ti S_{\l\i}(e(u))))$ 
\nl
or that $\tS_{\l\i}(\ph(e(u))=\ph(\tS_{\l\i}(e(u)))$. This follows from (c).

\mpb

For $u\in\fU_{w_I}(K)$ we show: 

(e) $S_\l(\Ph'(u))=\Ph'(S_{\l\i}(u))$.
\nl
This follows from (d), 4.1(e) and the identity $S_\l\ti\Ps=\ti\Ps S_\l$ (which follows from definitions).

\subhead 4.4\endsubhead
For $\l:I@>>>K,\l':I@>>>K$ we have $S_\l S_{\l'}=S_{\l\l'}$. Using this, together with 4.3(d),(e) and the 
equalities $\Ph^2=1$, $\Ph'{}^2=1$, we see that

(a) $(S_\l\Ph)^2=1$, $(S_\l\Ph')^2=1$
\nl
as maps $\fU_{w_I}(K)@>>>\fU_{w_I}(K)$.

\head 5. The positive part of a partial flag manifold\endhead
\subhead 5.1\endsubhead
We assume that $K,\kk$ are as in 1.1(i). 
Let $J\sub I$. Let $P_J^+$ be the subgroup of $G$ generated by $x_i(a)$ with $i\in I,a\in\kk$, by $y_i(a)$ with 
$i\in J,a\in\kk$ and by $T$. Let $P^-_J$ be the subgroup of $G$ generated by $y_i(a)$ with $i\in I,a\in\kk$, 
by $x_i(a)$ with $i\in J,a\in\kk$ and by $T$. Let $\cp_J$ be the set of subgroups of $G$ of the form $gP^-_Jg\i$ 
for some $g\in G$ or equivalently of the form $gP^+_{J^!}g\i$ for some $g\in G$. Here $J^!\sub I$ is the image of 
$J$ under the involution $i\m i^!$ of $I$ (see 2.1). Note that $P^+_\emp=B^+$, 
$P^-_\emp=B^-$, $\cp_\emp=\cb$. Let $\p_J:\cb@>>>\cp_J$ be the map which to any $B\in\cb$ associates the unique 
subgroup $P\in\cp_J$ such that $B\sub P$. We have $\p_J(B^-)=P_J^-$, $\p_J(B^+)=P_{J^!}^+$. Following \cite{L98}
we define $\cp(K)_{J,>0}=\p_J(\cb(K)_{>0})$, a subset of $\cp_J$.

Let $w=w_Iw_J$, $m=|w|=|w_I|-|w_J|$. Let $\ii=(i_1,i_2,\do,i_m)\in\co_w$. Define 
$f_\ii:K^m@>>>\cp_J$ by $(a_1,a_2,\do,a_m)\m uP_J^-u\i$ where 
$$u=x_{i_1}(a_1)x_{i_2}(a_2)\do x_{i_m}(a_m)\in U^+.$$ 
Define $f'_\ii:K^m@>>>\cp_J$ by $(a_1,a_2,\do,a_m)\m u'P_{J^!}^+u'{}\i$ where 
$$u'=y_{i_m}(a_m)y_{i_{m-1}}(a_{m-1})\do y_{i_1}(a_1)\in U^-.$$ We show:

(a) {\it $f_\ii,f'_\ii$ are injective.}
\nl
Assume that $(a_1,a_2,\do,a_m)\in K^m,(a'_1,a'_2,\do,a'_m)\in K^m$ are such that, 
setting 
$$u=x_{i_1}(a_1)x_{i_2}(a_2)\do x_{i_m}(a_m)\in U^+, u'=x_{i_1}(a'_1)x_{i_2}(a'_2)\do x_{i_m}(a'_m)\in U^+,$$
 we have $uP_J^-u\i=u'P_J^-u'{}\i$ that is, $u=u'g$ for some $g\in P_J^-$. 
We have $u\in G_w,u'\in G_w,g\in G_y$ (see 2.1) for some $y\in W_J$. Since $|wy|=|w|+|y|$, we have
$u'g\in G_{wy}$. Since $u=u'g$, it follows that $w=wy$ hence $y=1$. Thus, $g\in B^-$.
We have $u\in U^+,u'\in U^+$ hence $g=u'{}\i u\in U^+\cap B^-=\{1\}$. Thus, $u'=u$. From 
$x_{i_1}(a_1)x_{i_2}(a_2)\do x_{i_m}(a_m)=x_{i_1}(a'_1)x_{i_2}(a'_2)\do x_{i_m}(a'_m)$ 
we deduce that $(a_1,a_2,\do,a_m)=(a'_1,a'_2,\do,a'_m)$ (see 2.3(a)). This proves (a) for $f_\ii$.

Next we assume that $(a_1,a_2,\do,a_m)\in K^m,(a'_1,a'_2,\do,a'_m)\in K^m$ are such that, setting 
$$u=y_{i_m}(a_m)\do y_{i_2}(a_2)y_{i_1}(a_1)\in U^-,u'=y_{i_m}(a'_m)\do y_{i_2}(a'_2)y_{i_1}(a'_1)\in U^-,$$
we have $uP_{J^!}^+u\i=u'P_{J^!}^+u'{}\i$. Applying the involution $\io:G@>>>G$ (see 2.1), we deduce 
$\io(u)P_{J^!}^-\io(u\i)=\io(u')P_{J^!}^-\io(u'{}\i)$ where 
$$\io(u)=x_{i_m}(a_m)\do x_{i_2}(a_2)x_{i_1}(a_1)\in U^+,
\io(u')=x_{i_m}(a'_m)\do x_{i_2}(a'_2)x_{i_1}(a'_1)\in U^+.$$
It follows that $u_1P_J^-u_1\i=u'_1P_J^-u'_1{}\i$ where 
$$u_1=x_{i_m^!}(a_m)\do x_{i_2^!}(a_2)x_{i_1^!}(a_1)\in U^+,u'_1=x_{i_m^!}(a'_m)\do x_{i_2^!}(a'_2)x_{i_1^!}(a'_1)\in U^+.$$
Applying the first part of the argument to $(i_m^!,i_{m-1}^!,\do,i_1^!)$ instead of $\ii$
we deduce that $u_1=u'_1$ and $(a_1,a_2,\do,a_m)=(a'_1,a'_2,\do,a'_m)$. This proves (a) for $f'_\ii$.

We show:

(b) {\it The image of $f_\ii$ is equal to $\cp(K)_{J,>0}$; in particular it is independent of $\ii$.
The image of $f'_\ii$ is equal to $\cp(K)_{J,>0}$; in particular it is independent of $\ii$.}
\nl
Let $\jj=(i_1,i_2,\do,i_m,j_1,j_2,\do,j_s)$ where $(j_1,j_2,\do,j_s)\in\co_{w_J}$. We have $\jj\in\co_{w_I}$ 
hence $\cp(K)_{J,>0}$ is the set of all $gP_J^-g\i$ where 
$$g=x_{i_1}(a_1)x_{i_2}(a_2)\do x_{i_m}(a_m)x_{j_1}(b_1)x_{j_2}(b_2)\do x_{j_s}(b_s)$$
with $a_1,a_2,\do,a_m,b_1,b_2,\do,b_s$ in $K$. Since $x_{j_1}(b_1),x_{j_2}(b_2),\do,x_{j_s}(b_s)$ 
are in $P_J^-$, it follows that
$\cp(K)_{J,>0}$ is the set of all $gP_J^-g\i$ where $$g=x_{i_1}(a_1)x_{i_2}(a_2)\do x_{i_m}(a_m)$$
with $a_1,a_2,\do,a_m$ in $K$. Thus $\cp(K)_{J,>0}$ is the image of $f_\ii$ (see 4.1(a)). 
This proves the first sentence in (b).

Let $\jj'=(i_m,i_{m-1},\do,i_1,j_1^!,j_2^!,\do,j_s^!)$. We have $\jj'\in\co_{w_I}$ hence $\cp(K)_{J,>0}$
is the set of all $gP_{J^!}^+g\i$ where
$$g=y_{i_m}(a_m)\do y_{i_2}(a_2)y_{i_1}(a_1)y_{j^!_1}(b_1)y_{j^!_2}(b_2)\do y_{j^!_s}(b_s)$$
with $a_1,a_2,\do,a_m,b_1,b_2,\do,b_s$ in $K$.
Since $y_{j^!_1}(b_1),y_{j^!_2}(b_2),\do,y_{j^!_s}(b_s)$ are in $P_{J^!}^+$, it follows that
$\cp(K)_{J,>0}$ is the set of all $gP_{J^!}^+g\i$ where
$$g=y_{i_m}(a_m)\do y_{i_2}(a_2)y_{i_1}(a_1)$$ 
with $a_1,a_2,\do,a_m$ in $K$.
Thus, $\cp(K)_{J,>0}$ is the image of $f'_\ii$ (see 4.1(a)). This proves the second sentence in (b).

\subhead 5.2\endsubhead
Let $\ii,\jj,\jj',j_1,j_2,\do,j_s$ be as in 5.1(b) and its proof. Define $\tf_\jj:K^\nu@>>>\cb(K)_{>0}$ by 
$$\tf_\jj(a_1,a_2,\do,a_m,b_1,b_2,\do,b_s)=gB^-g\i$$
where $g=x_{i_1}(a_1)x_{i_2}(a_2)\do x_{i_m}(a_m)x_{j_1}(b_1)x_{j_2}(b_2)\do x_{j_s}(b_s)$.
Define $\tf'_{\jj'}:K^\nu@>>>\cb(K)_{>0}$ by 
$$\tf'_{\jj'}(a_1,a_2,\do,a_m,b_1,b_2,\do,b_s)=g'B^+g'{}\i$$
where $g'=y_{i_m}(a_m)\do y_{i_2}(a_2)y_{i_1}(a_1)y_{j^!_1}(b_1)y_{j^!_2}(b_2)\do y_{j^!_s}(b_s)$.
Let $A:K^\nu@>>>K^\nu$ be the bijection given by $\tf'_{\jj'}A=\tf_\jj$. Let 
$\bar A:K^m@>>>K^m$ be the bijection given by $f'_\ii\bar A=f_\ii$. We show:

(a) {\it there is a commutative diagram
$$\CD
K^\nu@>A>> K^\nu\\
@VVV            @VVV\\
K^m@>\bar A>>K^m\\
\endCD$$
where the vertical maps are given by }
$$(a_1,a_2,\do,a_m,b_1,b_2,\do,b_s)\m(a_1,a_2,\do,a_m).$$
Let $(a_1,a_2,\do,a_m,b_1,b_2,\do,b_s)\in K^\nu$ and let
$$(a'_1,a'_2,\do,a'_m,b'_1,b'_2,\do,b'_s)\in K^\nu$$ be its image under $A$. 
We must show that $(a'_1,a'_2,\do,a'_m)\in K^m$ is the image of $(a_1,a_2,\do,a_m)\in K^m$ under 
$\bar A$.
Let $$g=x_{i_1}(a_1)x_{i_2}(a_2)\do x_{i_m}(a_m)x_{j_1}(b_1)x_{j_2}(b_2)\do x_{j_s}(b_s),$$
$$g'=y_{i_m}(a'_m)\do y_{i_2}(a'_2)y_{i_1}(a'_1)y_{j^!_1}(b'_1)y_{j^!_2}(b'_2)\do y_{j^!_s}(b'_s).$$
 We have $gB^-g\i=g'B^+g'{}\i$ that is $g_1\tB g_1\i=g'_1\tB'g'_1{}\i$, $\tB=g_2B^-g_2\i$, 
$\tB'=g'_2B^+g'_2{}\i$, 
where $$g_1=x_{i_1}(a_1)x_{i_2}(a_2)\do x_{i_m}(a_m), g'_1=y_{i_m}(a'_m)\do y_{i_2}(a'_2)y_{i_1}(a'_1),$$
$$g_2=x_{j_1}(b_1)x_{j_2}(b_2)\do x_{j_s}(b_s),g'_2=y_{j^!_1}(b'_1)y_{j^!_2}(b'_2)\do y_{j^!_s}(b'_s).$$
Since $x_{j_1}(b_1),x_{j_2}(b_2),\do,x_{j_s}(b_s)$ are contained in $P_J^-$ and
$$y_{j^!_1}(b'_1),y_{j^!_2}(b'_2),\do, y_{j^!_s}(b'_s)$$ are contained in $P_{J^!}^+$
we see that $\tB\sub P_J^-$ and $\tB'\sub P_{J^!}^+$.
Hence 
$$g_1\tB g_1\i\sub g_1P_J^-g_1\i\text{ and }g'_1\tB'g'_1{}\i\sub g'_1P_{J^!}^+g'_1{}\i.$$
Since a subgroup in $\cb$ is contained in a unique subgroup in $\cp_J$, we deduce that
$g_1P_J^-g_1\i=g'_1P_{J^!}^+g'_1{}\i$. This proves (a).

Using (a) and the fact that $A:K^\nu@>>> K^\nu$ is admissible (see 4.1(a)), 
we see that $\bar A:K^m@>>>K^m$ is admissible. Now (a) remains a commutative diagram if $A$ is 
replaced by $A\i$ and $\bar A$ is replaced by $\bar A\i$. Using this and the 
fact that $A\i:K^\nu@>>> K^\nu$ is admissible (see 4.1(a)), we see that $\bar A\i:K^m@>>>K^m$ is admissible.

\subhead 5.3\endsubhead
Let $\ii,f_\ii:K^m@>>>\cp(K)_{J,>0}$ be as in 5.1 and let $\jj,\tf_\jj:K^\nu@>>>\cb(K)_{>0}$ be as in 5.2.
Let $\ii^1$ be another sequence in $\co_w$; let $\ff_{\ii^1}, \jj^1,\tf_{\jj^1}$ be defined in terms of
$\ii^1$ in the same way as $\ff_\ii, \jj,\tf_\jj$ were defined in terms of $\ii$. 
Let $A':K^\nu@>>>K^\nu$ be the bijection given by $\tf_{\jj^1}A'=\tf_\jj$. Let 
$\bar A':K^m@>>>K^m$ be the bijection given by $f_{\ii^1}\bar A'=f_\ii$. 
From the definitions we have a commutative diagram
$$\CD
K^\nu@>A'>> K^\nu\\
@VVV            @VVV\\
K^m@>\bar A'>>K^m\\
\endCD$$
where the vertical maps are as in 5.2(a). Since $A'$ and $A'{}\i$ are admissible it follows that
$\bar A'$ and $\bar A'{}\i$ are admissible. This, together with the results in 5.2 implies that

(a) {\it the bijections $f_\ii:K^m@>>>\cp(K)_{J,>0}$, $f'_\ii:K^m@>>>\cp(K)_{J,>0}$ (with $\ii\in\co_w$) define a 
positive $K$-structure on $\cp(K)_{J,>0}$.}

\mpb

From 3.4(b) and the definitions we deduce:

(b) {\it There is a well defined action of the monoid $G(K)$ on $\cp(K)_{J,>0}$ given by 
$g_1:gP_J^-g\i\m g_1gP_J^-g\i g_1\i$ where $gP^-_Jg\i\in\cp(K)_{J,>0}$.}

\mpb

From the definitions we see that for $w,w'$ in $W$, the map 
$$G_{w,-w'}(K)\T\cp(K)_{J,>0}@>>>\cp(K)_{J,>0}$$ (restriction of the $G(K)$-action
in (b)) is admissible. Hence the map 

(c) $G(K)\T\cp(K)_{J,>0}@>>>\cp(K)_{J,>0}$ 
\nl
given by the $G(K)$-action is admissible. 

Now let $K'=\RR(t)_{>0}$. We set $\cp(\cz)_{J,>0}=\un{\cp(K')_{J,>0}}$. By passage to zones in (c) we obtain
an admissible map $G(\cz)\T\cp(\cz)_{J,>0}@>>>\cp(\cz)_{J,>0}$ which is an action of $G(\cz)$.

\head 6. Coordinate rings\endhead
\subhead 6.1\endsubhead
In this section we assume that we are in the setup of 2.1 with $\kk$ an algebraically closed field of any 
characteristic. For any irreducible algebraic variety $V$ over $\kk$ be denote by $O(V)$ the algebra of regular 
functions $V@>>>\kk$ and by $[O(V)]$ the quotient field of $O(V)$.

For any $\ii=(i_1,i_2,\do,i_\nu)\in\co_{w_I}$, the map $f_\ii:\kk^\nu@>>>U^+$ given by
$$(a_1,a_2,\do,a_\nu)\m x_{i_1}(a_1)x_{i_2}(a_2)\do x_{i_\nu}(a_\nu)$$ defines a field isomorphism
$f_\ii^*:[O(U^+)]@>>>[O(\kk^\nu)]$.
We fix a numbering $I=\{1,2,\do,k,k+1,\do,r\}$ where $s_i,s_j$ commute if $i\le k\ge j$ or if $i>k<j$. Let 
$\jj'=(1,2,\do,k)$, $\jj''=(k+1,k+2,\do,r)$. Let $\jj\in\co_{w_I}$ be the concatenation $\jj'\jj''\jj'\jj''\do$
(of length $\nu$) and let $\jj'\in\co_{w_I}$ be the concatenation $\jj''\jj'\jj''\jj'\do$
(of length $\nu$). We state the following conjecture.

(a) {\it Let $F\in[O(U^+)]$. We have $F\in O(U^+)$ if and only if $f_\jj^*(F)\in O(\kk^\nu)$ and 
$f_{\jj'}^*(F)\in O(\kk^\nu)$.}
\nl
A consequence of (a) is that the algebra $O(U^+)$ can be identified with the algebra consisting of all
$(h,h')\in O(\kk^\nu)\T O(\kk^\nu)$ such that under the field isomorphism $[O(\kk^\nu)]@>\si>>[O(\kk^\nu)]$ 
induced by the birational isomorphism $f_\jj f_{\jj'}\i$, $h,h'$ correspond to each other. Thus the algebra
$O(U^+)$ is described completely in terms of the transition function between two charts in
the positive $\RR_{>0}$-structure attached to $U^+$.

We will verify (a) in two examples.

Assume first that $G$ is of type $A_2$. We can assume that $I=\{1,2\}$, $\jj=(1,2,1)$, $\jj'=(2,1,2)$.
We identify $\kk^3$ with $U^+$ as varieties by $(A,B,C)\m x_1(A)x_{12}(B)x_2(C)$
where $x_{12}:\kk@>>>U^+$ is a root subgroup corresponding to the nonsimple positive root.
In these coordinates the map $f_\jj$ (resp. $f_{\jj'}$) is given by $(a,b,c)\m(A,B,C)$ where
$A=a+c,B=bc,C=b$ (resp.  $A=b,B=ab,C=a+c$). The inverse (rational) map to $f_\jj$ is 
$(A,B,C)\m((AC-B)/C,C,A/C)$; it is regular on the complement of the hypersurface $C=0$. The inverse 
(rational) map to $f_{\jj'}$ is $(A,B,C)\m(B/A,A,(AC-B)/A)$; it is regular on the complement of the 
hypersurface $A=0$. If $F\in[O(U^+)]$ is regular after taking inverse image under $f_\jj$ and under $f_{\jj'}$ 
then it is regular outside the set $\{(A,B,C);A=C=0\}$. Since this set has codimension $2$ in $U^+$, $F$ must be
 regular. Thus (a) holds in this case. 

Next we assume that $G$ is of type $A_3$.
We can assume that $I=\{1,3,2\}$, $\ii^1=(1,3,2,1,3,2)$, $\ii^2=(2,1,3,2,1,3)$.
We identify $\kk^6$ with $U^+$ as varieties by
$$(X,Y,Z,U,V,W)\m x_1(X)x_3(Y)x_{213}(Z)x_{21}(U)x_{23}(V)x_2(W)$$
where $x_{21},x_{23},x_{123}$ are root subgroups $\kk@>>>U^+$ corresponding to the nonsimple positive roots.
In these coordinates the map $f_\jj$ (resp. $f_{\jj'}$) is given by 
$$(a,b,c,d,e,f)\m(X,Y,Z,U,V,W)$$ 
where
$$X=a+d,Y=b+c,Z=cde,U=cd,V=ce,W=c+f$$
 (resp. 
$$\align&X=c+e,Y=b+f,Z=def+a(b+f)(c+e),U=a(c+e)+de,\\&V=ab+af+df,W=a+d).\endalign$$
The inverse (rational) map to $f_\jj$ is 
$$(X,Y,Z,U,V,W)\m((XV-Z)/V,(UX-Z)/U,VU/Z,Z/V,Z/U,(WZ-VU)/Z);$$
it is regular on the complement of the hypersurface $XUV=0$.
The inverse (rational) map to $f_{\jj'}$ is $(X,Y,Z,U,V,W)\m(a,b,c,d,e,f)$ where
$$a=(ZW-UV)/(XYW-XV-UY+Z),$$
$$b=(XYW-XV-UY+Z)/(XW-U),$$
$$c=(XYW-XV-UY+Z)/(WY-V),$$
$$d=(XW-U)(WY-V)/(XYW-XV-UY+Z),$$ 
$$e=(UY-Z)/(WY-V),\qua f=(XV-Z)/(XW-U);$$
thisis regular on the complement of the hypersurface $(WY-V)(XYW-XV-UY+Z)(XW-U)=0$.
To prove (a) in this case it is then enough to observe that the last two hypersurfaces have
interesection of codimension $\ge2$ (they don't have an irreducible component in common).

\subhead 6.2\endsubhead
We write $i^1,i^2,\do,i^r$ for the order $1,2,\do,r$ on $I$.
Let $\ii=(i_1,i_2,\do,i_\nu)\in\co_{w_I}$,
$\ii'=(i'_1,i'_2,\do,i'_\nu)\in\co_{w_I}$. We consider the map $f_{\ii,\ii'}:\kk^\nu\T(\kk^*)^r\T\kk^\nu@>>>G$ 
given by
$$\align&(a_1,a_2,\do,a_\nu,b_1,b_2,\do,b_r,c_1,c_2,\do,c_\nu)\m\\& x_{i_1}(a_1)x_{i_2}(a_2)\do x_{i_\nu}(a_\nu)
i^1(b_1)i^2(b_2)\do i^r(b_r)y_{i'_1}(c_1)y_{i'_2}(c_2)\do y_{i'_\nu}(c_\nu)\endalign$$
and the map $\tf_{\ii,\ii'}:\kk^\nu\T(\kk^*)^r\T\kk^\nu@>>>G$ given by
$$\align&(a_1,a_2,\do,a_\nu,b_1,b_2,\do,b_r,c_1,c_2,\do,c_\nu)\m \\&y_{i_1}(a_1)y_{i_2}(a_2)\do y_{i_\nu}(a_\nu)
i^1(b_1)i^2(b_2)\do i^r(b_r)x_{i'_1}(c_1)x_{i'_2}(c_2)\do x_{i'_\nu}(c_\nu).\endalign$$
These defines field isomorphisms $f_{\ii,\ii'}^*,\tf_{\ii,\ii'}^*$ from $[O(G)]$ to 
$[O(\kk^\nu\T(\kk^*)^r\T\kk^\nu)]$. 
We state the following conjecture.

(a) {\it Let $\jj,\jj'$ be as in 6.1. Let $F\in[O(G)]$. We have $F\in O(G)$ if and only if 
$f_{\jj,\jj'}^*(F)\in O(\kk^\nu\T(\kk^*)^r\T\kk^\nu)$ and 
$\tf_{\jj',\jj}^*(F)\in O(\kk^\nu\T(\kk^*)^r\T\kk^\nu)$.}
\nl
More precisely, the inverse of the rational map $f_{\jj,\jj'}$ is regular on the complement of a hypersurface in $G$ and 
the inverse of the rational map $\tf_{\jj',\jj}$ is regular on the complement of another hypersurface in $G$
 and it should be true that
these two hypersurfaces have intersection of codimension $\ge2$ (they don't have a common irreducible
component). 

A consequence of (a) is that the algebra $O(G)$ can be identified with the algebra consisting of all
$(h,h')\in O(\kk^\nu\T(\kk^*)^r\T\kk^\nu)\T O(\kk^\nu\T(\kk^*)^r\T\kk^\nu)$
such that under the field isomorphism $[O(\kk^\nu\T(\kk^*)^r\T\kk^\nu)]@>\si>>[O(\kk^\nu\T(\kk^*)^r\T\kk^\nu)]$
induced by the birational isomorphism $f_{\jj,\jj'}\tf_{\jj',\jj}\i$, $h,h'$ correspond to each other. Thus the 
algebra $O(G)$ is described completely in terms of the transition function between two charts in
the positive $\RR_{>0}$-structure attached to $G$.

We verify the conjecture when $G$ is of type $A_1$. Let $V$ be the variety of all
$\left(\sm A&B\\C&D\esm\right)$ in $\kk^4$ such that $AD-BC=1$. 
Consider the map $s:\kk\T\kk^*\T\kk@>>>V$ given by
$(a,x,b)\m\left(\sm A&B\\C&D\esm\right)$ where $A=x+abx\i,B=ax\i,C=bx\i,D=x\i$ and the map
$s':\kk\T\kk^*\T\kk@>>>V$ given by $(a,x,b)\m\left(\sm A&B\\C&D\esm\right)$ where $A=x,B'=ax,C=bx,D=abx+x\i$.
The inverse of $s$ as a rational map is $\left(\sm A&B\\C&D\esm\right)\m(B/D,1/D,C/D)$.
The inverse of $s'$ as a rational map is $\left(\sm A&B\\C&D\esm\right)\m(B/A,A,C/A)$.
It remains to note that a rational map $V@>>>\kk$
which is regular on the subset $D\ne0$ and is regular on the subset $A\ne0$ is regular everywhere.

We verify the conjecture when $G$ is of type $A_2$. Let $V$ be the variety of all
$\left(\sm A&B&C\\E&D&F\\G&H&J\esm\right)$ in $\kk^9$ with determinant $1$.
Consider the map $s:\kk^3\T(\kk^*)^2\T\kk^3@>>>V$ given by
$$(a,b,c,\e,\d,a',b',c')\m
\left(\sm\e+(a+c)\e\i\d b'+ab\d\i a'b'&(a+c)\e\i\d +ab\d\i(a'+c')&ab\d\i\\ \e\i\d b'+b\d\i a'b'&\e\i\d+b\d\i(a'+c')
&b\d\i\\\d\i a'b'&\d\i(a'+c')&\d\i\esm\right).$$
Consider the map $s':\kk^3\T(\kk^*)^2\T\kk^3@>>>V$ given by
$$(a,b,c,\e,\d,a',b',c')\m
\left(\sm\e&\e(a'+c')&\e a'b'\\b\e&b\e(a'+c')+\e\i\d&b\e a'b'+\e\i\d b'\\ab\e&ab\e(a'+c')+(a+c)\e\i\d
&ab\e a'b'+(a+c)\e\i\d b'+\d\i\esm\right).$$
The inverse of $s$ as a rational map is 
$$\left(\sm A&B&C\\E&D&F\\G&H&J\esm\right)\m(a,b,c,\e,\d,a',b',c')$$
where
$$\align&a=\fra{C}{F},b=\fra{F}{J},c=\fra{(BF-CE)J}{(EJ-FH)F},\e=\fra{1}{EJ-FH}),\d=\fra{1}{J},\\&
a'=\fra{(EJ-FH)G}{(DJ-FG)J},b'=\fra{DJ-FG}{EJ-FH},c'=\fra{DH-GE}{DJ-FG}.\endalign$$
The inverse of $s'$ as a rational map is 
$$\left(\sm A&B&C\\E&D&F\\G&H&J\esm\right)\m(a,b,c,\e,\d,a',b',c')$$
where
$$\align&a=\fra{G}{D},b=\fra{D}{A},c=\fra{(HD-EG)A}{(AE-BD)D},\e=A,\d=AE-BD,\\&a'=\fra{(AE-BD)C}{(AF-CD)A},
b'=\fra{AF-CD}{AE-BD},c'=\fra{BF-CE}{AF-CD}.\endalign$$
It remains to note that a rational map $V@>>>\kk$ which is regular on the subset defined by 
$FJ(EJ-FH)(DJ-FG)\ne0$ and is regular on the subset defined by $AD(AE-BD)(AF-CD)\ne0$ is regular everywhere.

\subhead 6.3\endsubhead
Let $i^1,i^2,\do,i^r$ be as 6.2. Let $\ii=(i_1,i_2,\do,i_\nu)\in\co_{w_I}$.
We consider the map $f_\ii:\kk^\nu\T(\kk^*)^r@>>>G/U^-$ given by
$$\align&(a_1,a_2,\do,a_\nu,b_1,b_2,\do,b_r)\m\\& x_{i_1}(a_1)x_{i_2}(a_2)\do x_{i_\nu}(a_\nu)
i^1(b_1)i^2(b_2)\do i^r(b_r)U^-\endalign$$
and the map $\tf_\ii:\kk^\nu\T(\kk^*)^r@>>>G/U^-$ given by
$$\align&(a_1,a_2,\do,a_\nu,b_1,b_2,\do,b_r)\m \\&y_{i_1}(a_1)y_{i_2}(a_2)\do y_{i_\nu}(a_\nu)
i^1(b_1)i^2(b_2)\do i^r(b_r)\dw_IU^-.\endalign$$
These define field isomorphisms $f_\ii^*,\tf_\ii^*$ from $[O(G/U^-)]$ to $[O(\kk^\nu\T(\kk^*)^r]$. 

Let $\jj$ be as in 6.1. We state the following conjecture.

(a) {\it  Let $F\in[O(G/U^-)]$. We have $F\in O(G/U^-)$ if and only if $f_\jj^*(F)\in O(\kk^\nu\T(\kk^*)^r)$ and 
$\tf_\jj^*(F)\in O(\kk^\nu\T(\kk^*)^r)$.}
\nl
More precisely, the inverse of the rational map $f_\jj$ is regular on the complement of a hypersurface in $G/U^-$
 and the inverse of the rational map $\tf_\jj$ is regular on the complement of another 
hypersurface in $G/U^-$ and it should be true that these two hypersurfaces have intersection of codimension 
$\ge2$ (they don't have a common irreducible component). 

A consequence of (a) is that the algebra $O(G/U^-)$ can be identified with the algebra consisting of all
$(h,h')\in O(\kk^\nu\T(\kk^*)^r)\T O(\kk^\nu\T(\kk^*)^r)$ such that $h,h'$ correspond to each other under the 
field isomorphism $\k^*:[O(\kk^\nu\T(\kk^*)^r)]@>\si>>[O(\kk^\nu\T(\kk^*)^r)]$ induced by the birational 
isomorphism $\k=f_\jj\tf_\jj\i$.  Thus the algebra $O(G/U^-)$ is described completely in terms of the transition 
function between two charts in the positive $\RR_{>0}$-structure attached to $G/U^-$.

We verify the conjecture when $G$ is of type $A_2$. Let $V$ be the variety of all
$\left(\sm A&B&C\\A'&B'&C'\esm\right)$ in $\kk^6$ such that $(A,B,C)\ne0,(A',B',C')\ne0$, $AA'+BB'+CC'=0$.
(We can identify $V=G/U^-$.)

Consider the map $s:\kk^3\T(\kk^*)^2@>>>V$ given by
$$(a,b,c,\e,\d)\m\left(\sm \d\i bc&\d\i (a+c)&\d\i    \\ \e\i&-e\i b&\e\i ab\esm\right)$$
and the map $s':\kk^3\T(\kk^*)^2@>>>V$ given by
$$(a,b,c,\e,\d)\m\left(\sm \e&b\e&ab\e\\ \d bc&-\d(a+c)&\d\esm\right).$$
The inverse of $s$ as a rational map is 
$$\left(\sm A&B&C\\A'&B'&C'\esm\right)\m(a,b,c,\e,\d)$$
where $a=-C'/A',b=-B'/A',c=-AA'/(CB'),\e=1/A',\d=1/C$.
The inverse of $s'$ as a rational map is 
$$\left(\sm A&B&C\\A'&B'&C'\esm\right)\m(a,b,c,\e,\d)$$
where $a=C/B$, $b=B/A$, $c=A'A/(C'B)$, $\e=A$, $\d=C'$.
It remains to note that a rational map $V@>>>\kk$ which is regular on the subset defined by 
$A'B'C\ne0$ and is regular on the subset defined by $ABC'\ne0$ is regular everywhere.

\subhead 6.4\endsubhead
We preserve the setup of 6.3. Now $T$ acts on $G/U^-$ by $t:gU^-\m gtU^-$. For any $\l\in\cx$ we set 
$$O(G/U^-)_\l=\{f\in O(G/U^-);f(gt\i U^-)=\l(t)f(gU^-)\qua\frl g\in G,t\in T\}.$$
 It is known that $O(G/U^-)_\l$ is a 
finite dimensional $\kk$-vector space and it is $\ne0$ if and only if $\l\in\cx^+$. Moreover as a $\kk$-vector 
space we have $O(G/U^-)=\op_{\l\in\cx^+}O(G/U^-)_\l$. The left translation by $G$ induces a $G$-action on 
$O(G/U^-)$ and this keeps stable each of the subspaces $O(G/U^-)_\l$ with $\l\in\cx^+$ which (in the case where
the characteristic of $\kk$ is zero) becomes an irreducible representation of $G$ (Borel-Weil).

Consider the $T$-action on $\kk^\nu\T(\kk^*)^r$ given by
$$t:(a_1,a_2,\do,a_\nu,b_1,b_2,\do,b_r)\m(a_1,a_2,\do,a_\nu,b_1t_1,b_2t_2,\do,b_rt_r)$$
where $t=i^1(t_1)\do i^r(t_r)$, $(t_1,\do,t_r)\in K^r$. This induces a $T$-action on $O(\kk^\nu\T(\kk^*)^r)$; 
for $\l\in\cx$ we denote by $O(\kk^\nu\T(\kk^*)^r)_\l$ the subspace of $O(\kk^\nu\T(\kk^*)^r)$ on which $T$ acts 
according to $\l$. 
We have an isomorphism $e_\l:O(\kk^\nu\T(\kk^*)^r)_\l@>\si>>O(\kk^\nu)$ given by $f\m f'$ where 
$$f'(a_1,a_2,\do,a_\nu)=f(a_1,a_2,\do,a_\nu,1,1,\do,1).$$
Now $f_\jj:\kk^\nu\T(\kk^*)^r@>>>G/U^-$ is $T$-equivariant hence $f_\jj^*$ defines an 
(injective) linear map $f_{\jj,\l}^*:O(G/U)_\l@>>>O(\kk^\nu\T(\kk^*)^r)_\l$. Similarly, 
$\tf_\jj:\kk^\nu\T(\kk^*)^r)_\l@>>>G/U^-$ is $T$-equivariant for the $T$-action on $\kk^\nu\T(\kk^*)^r$ 
described above and the $T$-action on $G/U^-$ given by $gU^-\m gw_ItW_i\i U^-$ hence $\tf_\jj^*$ defines an 
(injective) linear map $\tf_{\jj,\l}^*:O(G/U)_\l@>>>O(\kk^\nu\T(\kk^*)^r)_{-\l^!}$ with $\l^!$ as in 4.2(b).
For $\l\in\cx^+$, the following statement is a consequence of the conjecture 6.3(a):

(a) {\it  Let $F\in[O(G/U^-)]$. We have $F\in O(G/U^-)_\l$ if and only if 
$f_\jj^*(F)\in O(\kk^\nu\T(\kk^*)^r)_\l$ and $\tf_\jj^*(F)\in O(\kk^\nu\T(\kk^*)^r)_{-\l^!}$.}
\nl
As in 6.3, this implies:

(b) {\it $O(G/U^-)_\l$ can be identified with the vector space consisting of all
$(h,h')\in O(\kk^\nu\T(\kk^*)^r)_\l\T O(\kk^\nu\T(\kk^*)^r)_{-\l^!}$ such that $h,h'$ correspond to each other 
under the field isomorphism $\k^*$ in 6.3.}
\nl
Equivalently, 

(c) {\it $O(G/U^-)_\l$ can be identified with the vector space consisting of all
$(h_0,h'_0)\in O(\kk^\nu)\T O(\kk^\nu)$ such that $e_\l\i(h),e_{-\l^!}\i(h')$ correspond to each other 
under the field isomorphism $\k^*$ in 6.3.}
\nl
Thus, the vector space $O(G/U^-)_\l$ is described completely in terms of the transition 
function between two charts in the positive $\RR_{>0}$-structure attached to $G/U^-$.

We consider for example the case where $G$ is of type $A_2$ and $I=\{i^1,i^2\}$. Let 
$$(a,b,c,\e,\d),(a',b',c',\e',\d')$$
in $\kk^3\T(\kk^*)^2$ be such that $s(a,b,c,\e,\d)=s'(a',b',c',\e',\d'),$ (notation of 6.3), that is
$$\left(\sm \d\i bc&\d\i (a+c)&\d\i\\ \e\i&-e\i b&\e\i ab\esm\right)=
\left(\sm\e'&b'\e'&a'b'\e'\\ \d'b'c'&-\d'(a'+c')&\d'\esm\right).$$
We have
$$a'=1/(a+c),b'=(a+c)/(bc),c'=c/(a(a+c)),\e'=\d\i bc,\d'=\e\i ab.$$
Assuming that $\la i^1,\l\ra=x\in\NN,\la i^2,\l\ra=y\in\NN$, we see that $O(G/U^-)_\l$ is identified with the 
vector space $\cv$ of polynomials $\sum_{i,j,k\text{ in }\NN}n_{i,j,k}a^ib^jc^k$ in $a,b,c$ (with 
$n_{i,j,k}\in\kk$) such that
$$\sum_{i,j,k\text{ in }\NN}n_{i,j,k}(a+c)^{-i}(a+c)^jb^{-j}c^{-j}(a+c)^{-k}c^ka^{-k}(bc)^x(ab)^y\tag a$$
is a polynomial in $a,b,c$. If for example $x=1,y=0$, then by direct computation we see that $\cv$ has a 
basis consisting of the polynomials $1,b,ab$. (In this case $O(G/U^-)_\l$ is known to have dimension $3$, 
which agrees with our computation.) If $x=1,y=1$, then by direct computation we see that $\cv$ has a basis 
consisting of the polynomials $1,b,a+c,ab,bc,b^2c,ab(a+c),ab^2c$. (In this case $O(G/U^-)_\l$ is known to have 
dimension $8$ which agrees with our computation.)

\head 7. Arnold's problem\endhead
\subhead 7.1\endsubhead
We assume that $K,\kk$ are as in 1.1(i) with $\kk=\RR,K=\RR_{>0}$. Let $K'=\RR(t)_{>0}$.
Let $\cb^*$ be the set of all $B\in\cb$ such that $(B^+,B)$, $(B,B^-)$, $(B^-,B^+)$ are in the same $G$-orbit 
under simultaneous conjugation. Now $\cb^*$ is naturally a real algebraic manifold. Let $[\cb^*]$ be the set of 
connected components of $\cb^*$. Arnold's problem asks to compute the number of elements of $[\cb^*]$. For 
$G=SL_n$, this was solved in \cite{SSV}: the number is $2,6,20,52$ if $n=2,3,4,5$ and 
is $2\T3^{n-1}$ if $n\ge6$. (For $n=2$, $\cb^*$ is a circle minus two points; thus it has two connected 
components.) For a general $G$, Arnold's problem was solved by Rietsch \cite{R97}, who found a combinatorial 
parametrization of the set $[\cb^*]$. We reformulate her parametrization by defining a map 
$\s:U^+_{w_I}(\cz)/2\to[\cb^*]$ (see 2.7) which is a bijection. To define $\s$ it is enough to define a map 
$\ti\s:U^+_{w_I}(\cz)\to[\cb^*]$ which is constant on the equivalence classes for $\si_2$ (see 2.7). Since 
$U^+_{w_I}(\cz)$ is the set of set of zones in $U^+_{w_I}(K')$, it is enough to associate to 
each zone an element of $[\cb^*]$. Let $x\in U^+_{w_I}(K')$. Let $\ii=(i_1,i_2,\do,i_\nu)\in\co_{w_I}$. We have 

(a) $x=x_{i_1}(a_1)x_{i_2}(a_2)\do x_{i_\nu}(a_\nu)$ 
\nl
where $a_s:t\m t^{e_s}f_s(t)/f'_s(t)$ and $f_s\in\RR[t],f'_s\in\RR[t]$ have constant term in $\RR_{>0}$.
We can find $\d\in\RR_{>0}$ such that $f_s(\d')>0$, $f'_s(\d')>0$ for any $s$ and any $\d'\in(-\d,0)$. 
Hence for $\d'\in(-\d,0)$, $a_s(\d')$ is well defined, is in $(-1)^{e_s}\RR_{>0}$ and
$$x(\d')=x_{i_1}(a_1(\d'))x_{i_2}(a_2(\d'))\do x_{i_\nu}(a_\nu(\d'))\in U^+$$
is well defined. For $\d'\in(-\d,0)$ we have $x(\d')B^-x(\d')\i\in\cb^*$; this is contained in a 
connected component $\fC(x)$ of $\cb^*$ which is independent of the choice of $\ii,\d,\d'$. 

\subhead 7.2\endsubhead
Now let $\tx\in U_{w_I}(K')$ be in the same zone as $x$ (as in 7.1(a)). We show that $\fC(\tx)=\fC(x)$.
We have $\tx=x_{i_1}(\ta_1)x_{i_2}(\ta_2)\do x_{i_\nu}(\ta_\nu)$ with $\ii$ as above where $\ta_s:\RR@>>>\RR$
are rational functions of the form $t\m t^{e_s}\tf_s(t)/\tf'_s(t)$ and $\tf_s\in\RR[t],\tf'_s\in\RR[t]$ 
have constant term in $\RR_{>0}$. We can find $\d_1\in\RR_{>0}$ such that $(cf_s+(1-c)\tf_s)(\d')>0$, 
$(cf'_s+(1-c)\tf'_s)(\d')>0$ for any $s$, any $c\in[0,1]$ and any $\d'\in(-\d_1,0)$. For $c\in[0,1]$ we set 
$a_{c,s}=t^{e_s}(cf_s+(1-c)\tf_s)/(cf'_s+(1-c)\tf'_s)\in K'$ and 
$$x_c=x_{i_1}(a_{c,1})x_{i_2}(a_{c,2})\do x_{i_\nu}(a_{c,\nu})\in U_{w_I}(K').$$
 Then for $c\in[0,1]$, $\d'\in(-\d_1,0)$, 
$$x_c(\d')=x_{i_1}(a_{c,1}(\d'))x_{i_2}(a_{c,2}(\d'))\do x_{i_\nu}(a_{c,\nu}(\d'))\in U^+$$
is well defined. For $\d'\in(-\d_1,0)$ we have $x_c(\d')B^-x_c(\d')\i\in\cb^*$ and the connected component of
$\cb^*$ containig it is independent of $c,\d_1,\d'$. In particular for $\d'\in(-\d_1,0)$, $x(\d')B^-x)\d')\i$ 
and $\tx(\d')B^-\tx)\d')\i$ belong to the same connected component of $\cb^*$, so that $\fC(x)=\fC(\tx)$.

Next let $x'=x_{i_1}(a'_1)x_{i_2}(a'_2)\do x_{i_\nu}(a'_\nu)$ where $a'_s:\RR@>>>\RR$ are rational functions of 
the form $t\m t^{e_s+2n_s}f_s(t)/f'_s(t)$ where $f_s,f'_s,e_s$ are as above and $n_s\in\NN$. We show that 
$\fC(x')=\fC(x)$. We can find $\d_2\in\RR_{>0}$ such that $(c+(1-c)\d'{}^{2n_h})f_s(\d')>0$, $f'_s(\d')>0$
for any $s$, any $c\in[0,1]$ and any $\d'\in(-\d_2,0)$. For $c\in[0,1]$ we set 
$a'_{c,s}=t^{e_s}((c+(1-c)t^{2n_s})f_s)/f'_s\in K'$ and 
$$x'_c=x_{i_1}(a'_{c,1})x_{i_2}(a'_{c,2})\do x_{i_\nu}(a'_{c,\nu})\in U_{w_I}(K').$$
Then for $c\in[0,1],\d'\in(-\d_2,0)$,
$$x'_c(\d')=x_{i_1}(a'_{c,1}(\d'))x_{i_2}(a'_{c,2}(\d'))\do x_{i_\nu}(a'_{c,\nu}(\d'))\in U^+$$ 
is  well defined. For $\d'\in(-\d_2,0)$ we have $x'_c(\d')B^-x'_c(\d')\i\in\cb^*$ and the connected component of
$\cb^*$ containig it is independent of $c,\d_2,\d'$. In particular for $\d'\in(-\d_1,0)$, $x'(\d')B^-x'(\d)\i$ 
and $x(\d')B^-x(\d')\i$ belong to the same connected component of $\cb^*$, so that $\fC(x)=\fC(x'$.

We see that $x\m\fC(x)$ is a well defined map from the set of zones of $U_{w_I}(K')$ to the set of connected 
components of $\cb^*$ and this map is constant on the equivalence classes on the set of zones for the equivalence
relation $\sim_2$ in 2.7. 
The resulting map $U^+_{w_I}(\cz)/2\to[\cb^*]$ can be identified with a bijection defined in 
\cite{R97} hence is itself a bijection. 

\head 8. The sets $\fU^\l_{w_I}(\cn),\ti\fU^\l_{w_I}(\cn)$\endhead
\subhead 8.1\endsubhead
Let $K$ be as in 1.1(i)-(iii). For any $i\in I$ we define $z_i:\fU_{w_I}(K)@>>>K$ as follows. Let $x\in\fU_{w_I}(K)$.
We choose $\ii=(i_1,i_2,\do,i_{\nu-1})\in\co_{w_Is_i}$ and we set $z_{i,\ii}(x)=a_\nu$ where 
$x=i_1^{a_1}i_2^{a_2}\do i_{\nu-1}^{a_{\nu-1}}i^{a_\nu}$ and $(a_1,a_2,\do,a_\nu)\in K^\nu$ is uniquely determined by $x$ (note that 
$(i_1,i_2,\do,i_{\nu-1},i)\in\co_{w_I}$). We show that $z_{i,\ii}(x)$ is independent of the choice of $\ii$. Using 
Iwahori's lemma,
we see that it is enough to show that $z_{i,\ii}(x)=z_{i,\ii'}(x)$ if $\ii,\ii'$ in $\co_{w_Is_i}$ 
are connected by a single braid move. In this case the desired result follows from 2.9(ii). We set
$z_i(x)=z_{i,\ii}(x)$ where $\ii$ is any sequence in $\co_{w_Is_i}$.

\subhead 8.2\endsubhead
In the remainder of this section we assume that $K=\cz$.
Let $\l\in\cx^+$. The function $I@>>>\NN$, $i\m\la i,\l\ra$ is denoted again by $\l$.
Define $\l^!:I@>>>\NN$ by $\l^!(i)=\l(i^!)$ ($i^!$ as in 2.1). Let 
$$\fU^\l_{w_I}(\cn)=\{x\in\fU_{w_I}(\cn);z_i(x)\le\l(i)\text{ for any }i\in I\}.$$
It is known that $\fU^\l_{w_I}(\cn)$ is naturally an indexing set for the canonical basis of the simple
$G$-module $\L_\l$ (see 2.1) with $G$ over $\kk=\RR$. (This is proved in \cite{L90} assuming that $G$ is 
simply laced; see \cite{L92}, \cite{L11} for the reduction of the general case to the simply laced case.) 
According to \cite{L97, 4.9}, the bijection $S_\l\Ph:\fU_{w_I}(\cz)@>>>\fU_{w_I}(\cz)$ (see 4.1, 4.3)
restricts to a bijection $\fU^\l_{w_I}(\cn)@>>>\fU^{\l^!}_{w_I}(\cn)$. In particular we have 

(a) $\fU^\l_{w_I}(\cn)\sub\ti\fU^\l_{w_I}(\cn)$
\nl
where 
$$\ti\fU^\l_{w_I}(\cn)=\{x\in\fU_{w_I}(\cn);S_\l\Ph(x)\in\fU_{w_I}(\cn).$$ 
The following conjecture is suggested by the results in 6.4:

(b) {\it The inclusion (a) is an equality.}
\nl
We will verify this conjecture in the case where $G$ is of type $A_1$ or $A_2$.
Assume first that $G$ is of type $A_1$ and $I=\{1\}$. We have $\fU^\l_{w_I}(\cn)=\{1^a;a\in\NN,a\le\l(1)\}$.
In our case $\Ph$ is given by $1^a\m 1^{-a}$ for $a\in\ZZ$ and 
$S_\l(1^{a'})=1^{a'+\l(1)}$ for $a'\in\ZZ$. Hence $S_\l\Ph(1^a)=1^{-a+\l(1)}$ and
$$\ti\fU^\l_{w_I}(\cn)=\{1^a;a\in\NN^3,a\le\l(1)\}=\fU^\l_{w_I}(\cn).$$
Thus (b) holds. 

Next we assume that $G$ is of type $A_2$ and $I=\{1,2\}$. We have 
$$\align&\fU^\l_{w_I}(\cn)=\{2^a1^b2^c;(a,b,c)\in\NN^3,c\le\l(2),a+b-\min(a,c)\le\l(1)\}\\&
=\{2^a1^b2^c;(a,b,c)\in\NN^3,a\le c,c\le\l(2),b\le\l(1)\}\sqc\\&
\{2^a1^b2^c;(a,b,c)\in\NN^3,a>c,c\le\l(2),a+b-c\le\l(1)\}.\endalign$$
In our case $\Ph$ is given by $2^a1^b2^c\m 1^{c-a-b}2^{-c}1^{-b}$ for $(a,b,c)\in\ZZ^3$ and 
$$S_\l(1^{a'}2^{b'}1^{c'})=1^{a'+\l(1)}2^{b'+\l(2)}1^{c'+\l(1)}$$ for $(a',b',c')\in\ZZ^3$. 
Hence $$S_\l\Ph(2^a1^b2^c)=1^{c-a-b+\l(1)}2^{-c+\l(2)}1^{-b+\l(1)}$$ and
$$\ti\fU^\l_{w_I}(\cn)=\{2^a1^b2^c;(a,b,c)\in\NN^3,a+b-c\le\l(1),c\le\l(2),b\le\l(1)\}.$$
If $a>c$ we have $b\le a+b-c$ hence the condition $b\le\l(1)$ is a consequence of the condition $a+b-c\le\l(1)$;
if $a\le c$ we have $a+b-c\le b$ hence the condition $a+b-c\le\l(1)$ is a consequence of the condition $b\le\l(1)$. We see
that in our case we have $\ti\fU^\l_{w_I}(\cn)=\fU^\l_{w_I}(\cn)$ and (b) holds.

\head 9. On $G(K)$ and conjugacy classes\endhead
\subhead 9.1\endsubhead
We assume that $K,\kk$ are as in 1.1(i) with $\kk=\RR$, $K=\RR_{>0}$. We denote by $G(\CC)$ the group of points 
of $G$ over $\CC$.

Let $g\in G(K)$. We show:

(a) {\it If $V\in\cc$ (see 2.1), then all eigenvalues of $g:V@>>>V$ are in $K$.}
\nl
By \cite{L94, 4.4}, there exists a sequence $g(k)\in G_{w_I,-w_I}(K)$, $k=1,2,\do$ such that $g=\lim_{k\m\iy}g(k)$ (limit in $G$). 
If (a) holds when $g$ is replaced by any $g(k)$ with $k\ge1$ then it would follow that the eigenvalues of $g:V@>>>V$ are in 
$K\cup\{0\}$. Since these eigenvalues are $\ne0$ ($g$ is invertible in $G$) they must then be in $K$. Thus it is 
enough to prove (a) assuming that $g\in G_{w_I,-w_I}(K)$. In this case, by the Gantmacher-Krein theorem in type 
$A$ and by \cite{L94, 5.6} in the general case, $g$ is regular and semisimple
and by \cite{L94, 8.10}, we have $g=uu'tu\i$ where $u\in U^-_{w_I}(K)$, $u'\in U^+_{w_I}(K)$, $t\in T(K)$. Thus we can assume that 
in (a), $g$ is regular, semisimple and $g=u't$ with $u'\in U^+$, $t\in T(K)$. In this case $g$ is conjugate to $t$ by an element of $U^+$.
Hence we can assume that in (a) we have $g=t\in T(K)$. In this case the eigenvalues of $g:V@>>>V$ are of the form $\mu(t)$ where
$\mu\in\cx$. It remains to observe that $\mu:T@>>>\RR^*$ carries $T(K)$ into $K$. (We use that for $i\in I,a\in K$ we have
$\mu(i(a))=a^{\la i,\mu\ra}\in K$.)

\mpb

Let $g\in G(K)$. By \cite{L94, 8.11}, there exists $B\in\cb$ such that $g\in B$. Let $B(\CC)$ be the Borel 
subgroup of $G(\CC)$ for which $B$ is the group of real points. It follows that the semisimple part $g_s$ of $g$ 
is contained in $T'(\RR)$ where $T'$ is a maximal torus of $B(\CC)$, defined and split over $\RR$. We show:

(b) {\it $g_s$ is contained in the identity component of $T'(\RR)$.}
\nl
Let $\cx'$ be the group of homomorphisms of algebraic groups $T'@>>>\CC^*$. For any $\mu\in\cx'$, $\mu(t)$ is an 
eigenvalue of $g_s:\L_\l@>>>\L_\l$ for some $\l\in\cx^+$ hence it is in $K$ (by (a)). An element of $T'(\RR)$ 
such that $\mu(t)\in K$ for any $\mu\in\cx'$ is necessarily in the identity component of $T'(\RR)$. This proves 
(b).

We show:

(c) {\it In the setup of (b), the centralizer $Z(g_s)$ of $g_s$ in $G(\CC)$ is the Levi subgroup of some 
parabolic subgroup of $G(\CC)$ defined over $\RR$.}
\nl
Let $R$ be the set of all $\mu\in\cx'$ which are roots of $G(\CC)$ with respect to $T'$. Let $R_0$ be the set of 
all $\mu\in R$ such that $\mu(g_s)=1$. It is enough to show that $R_0$ is the intersection of $R$ with a 
$\QQ$-subspace of $\QQ\ot\cx'$ or equivalently that: 

(d) if $\a\in R$ and $\a^c=\prod_{\mu\in R_0}\mu^{n_\mu}$ for some integers $n_\mu$ and some $c\in\ZZ_{>0}$, 
then $\a\in R_0$.
\nl
In the setup of (d) we have $\a(g_s)^c=1$. From (b) we see that $\a(g_s)\in\RR_{>0}$ so that from
$\a(g_s)^c=1$ we can deduce that $\a(g_s)=1$. This proves (d) and hence (c).

\subhead 9.2\endsubhead
Let $g\in G(K)$. Let $Z(g_s)$ be as in 9.1(c). Let $g_u$ be the unipotent part of $g$, so that $g_u\in Z(g_s)$.
We conjecture that there exists a subgroup $H$ of $Z(g_s)$ which is a Levi subgroup of a parabolic subgroup of 
$Z(g_s)$ such that $g_u$ is a regular unipotent element of $H$; in particular, any unipotent element in $G(K)$
is a regular unipotent element in a Levi subgroup of some parabolic subgroup of $G(\CC)$. (This last statement 
is obvious if $G$ is of type $A$ and can be shown to be true if $G$ is of type $D$, using the description
of unipotent elements in $G(K)$ given in \cite{L94, 6.6}.)

\head 10. A partition of $\BB$\endhead
\subhead 10.1\endsubhead
Let $w\in W$, $m=|w|$. For any $\ii=(i_1,i_2,\do,i_m)\in\co_w$ and any $\aa=(a_1,a_2,\do,a_m)\in\NN^m$ we set 
$[\ii,\aa]=i_{j_1}^1i_{j_2}^1\do i_{j_k}^1$ (product in $\fU(\{1\})$) where $j_1<j_2<\do<j_k$ is the sequence 
consisting of all $j\in[1,m]$ such that $a_j=0$. We show:

(a) {\it Assume that $\ii,\ii'$ in $\co_w$, $\aa,\aa'$ in $\ZZ^m$ are such that $\ii,\ii'$ are equal except at 
the indices $l,l+1,l+2$ where they are respectively $i,j,i$ and $j,i,j$ with $i:j=-1$ and that $\aa,\aa'$ are 
equal except at the indices $l,l+1,l+2$ where they are respectively $a,b,c$ and $a',b',c'$ with 
$a'=b+c-\min(a,c),b'=\min(a,c),c'=a+b-\min(a,c)$. Assume also that $\aa\in\NN^m$. Then $\aa'\in\NN^m$ and 
$[\ii,\aa]=[\ii',\aa']$ (equality in $\fU(\{1\})$).}
\nl
The fact that $\aa'\in\NN^m$ is immediate. Let $j_1<j_2<\do<j_t$ be the sequence consisting of all $j\in[1,l-1]$ 
such that $a_j=0$ or equivalently $a'_j=0$; let $j_{t'}<j_{t'+1}<\do<j_k$ be the sequence consisting of all 
$j\in[l+3,m]$ such that $a_j=0$ or equivalently $a'_j=0$. Let $y=i_{j_1}^1i_{j_2}^1\do i_{j_t}^1\in U^+(\{1\})$,
$y'=i_{j_{t'}}^1i_{j_{t'+1}}^1\do i_{j_k}^1\in U^+(\{1\})$. We have
$[\ii,\aa]=y[(i,j,i),(a,b,c)]y'$, $[\ii',\aa']=y[(j,i,j),(a',b',c')]y'$. Thus, to prove (a) we can assume that
$\ii=(i,j,i)$, $\ii'=(j,i,j)$, $\aa=(a,b,c),\aa'=(a',b',c')$. We consider a number of cases.

(I) $a=b=c=0$ so that $a'=b'=c'=0$. Then $[\ii,\aa]=i^1j^1i^1=j^1i^1j^1=[\ii',\aa']$.

(II) $a=b=0,c>0$ so that $a'>0,b'=c'=0$. Then $[\ii,\aa]=i^1j^1=[\ii',\aa']$.

(III) $a>0,b=c=0$ so that $a'=b'=0,c'>0$. Then $[\ii,\aa]=j^1i^1=[\ii',\aa']$.

(IV) $a=0,b>0,c=0$ so that $a'>0,b'=0,c'>0$. Then $[\ii,\aa]=i^1i^1=i^1=[\ii',\aa']$.

(V) $a=0,b>0,c>0$ so that $a'>0,b'=0,c'>0$. Then $[\ii,\aa]=i^1=[\ii',\aa']$.

(VI) $a>0,b>0,c=0$ so that $a'>0,b'=0,c'>0$. Then $[\ii,\aa]=i^1=[\ii',\aa']$.

(VII) $c>a>0,b=0$ so that $a'>0,b'>0,c'=0$. Then $[\ii,\aa]=j^1=[\ii',\aa']$.

(VIII) $a>c>0,b=0$ so that $a'=0,b'>0,c'>0$. Then $[\ii,\aa]=j^1=[\ii',\aa']$. 

(IX)  $a=c>0,b=0$ so that $a'=0,b'>0,c'=0$. Then $[\ii,\aa]=j^1=j^1j^1=[\ii',\aa']$.

(X) $a>0,b>0,c>0$ so that $a'>0,b'>0,c'>0$. Then $[\ii,\aa]=1=[\ii',\aa']$.
\nl
In each case, (a) is proved.

We show:

(b) {\it Assume that $\ii,\ii'$ in $\co_w$, $\aa,\aa'$ in $\ZZ^m$ are such that $\ii,\ii'$ are equal except at 
the indices $l,l+1$ where they are respectively $i,j$ and $j,i$ with $i:j=0$ and that $\aa,\aa'$ are 
equal except at the indices $l,l+1$ where they are respectively $a,b$ and $a',b'$ with $a'=b,b'=a$. Assume also 
that $\aa\in\NN^m$. Then $\aa'\in\NN^m$ and $[\ii,\aa]=[\ii',\aa']$ (equality in $\fU(\{1\})$).}
\nl
As in the proof of (a) we can assume that $\ii=(i,j)$, $\ii'=(j,i)$, $\aa=(a,b),\aa'=(a',b')$. 
We consider a number of cases.

(I) $a=b=0$ so that $a'=b'=0$. Then $[\ii,\aa]=i^1j^1=j^1i^1=[\ii',\aa']$.

(II) $a=0,b>0$ so that $a'>0,b'=0$. Then $[\ii,\aa]=i^1=[\ii',\aa']$.

(III) $a>0,b=0$ so that $a'=0,b'>0$. Then $[\ii,\aa]=j^1=[\ii',\aa']$.

(IV) $a>0,b>0$ so that $a'>0,b'>0$. Then $[\ii,\aa]=1=[\ii',\aa']$.
\nl
In each case, (b) is proved.

We show:

(c) {\it Assume that $\ii,\ii'$ in $\co_w$, $\aa,\aa'$ in $\NN^m$ satisfy $e_\ii(\aa)=e_{\ii'}(\aa')$ 
(equality in $\fU(\cz)$, $e_\ii,e_{\ii'}$ as in 2.9). Then $[\ii,\aa]=[\ii',\aa']$ (equality in $\fU(\{1\})$).}
\nl
Assume first that $G$ is simply laced.
By the Iwahori-Matsumoto lemma we can find a sequence $\ii=\ii^0,\ii^1,\do,\ii^s=\ii'$ in $\co_w$ 
such that for any $u\in\{0,1,\do,s-1\}$, $\ii^u,\ii^{u+1}$ are like $\ii,\ii'$ in (a) or (b).
We define a sequence $\aa=\aa^0,\aa^1,\do,\aa^s=\aa'$ in $\ZZ^m$ such that for any $u\in\{0,1,\do,s-1\}$, 
$\ii^u,\aa^u$ and $\ii^{u+1},\aa^{u+1}$ are like $\ii,\aa$ and $\ii',\aa'$ in (a) or (b).
Note that each $\aa^t$ with $0\le t\le s$ is automatically in $\NN^m$ (this follows by induction on $t$ from 
(a),(b)). From the definition we have $e_{\ii^u}(\aa^u)=e_{\ii^{u+1}}(\aa^{u+1})$.
From (a),(b), we see that $[\ii^u,\aa^u]=[\ii^{u+1},\aa^{u+1}]$ for any $u\in\{0,1,\do,s-1\}$.
It follows that $e_\ii(\aa)=e_{\ii'}(\aa^s)$ and $[\ii,\aa]=[\ii',\aa^s]$.
Hence we have $e_{\ii'}(\aa^s)=e_{\ii'}(\aa')$. Since $e_{\ii'}:\ZZ^m@>>>\fU_w(\cz)$ is a bijection
it follows that $\aa^s=\aa'$, so that $[\ii,\aa]=[\ii',\aa']$ and (c) is proved. 
If $G$ is not simply laced the proof is similar or it can be reduced to the simply laced by descent as in
\cite{L94}.

\mpb

We show:

(d) {\it Let $\ii=(i_1,i_2,\do,i_m)\in\co_w$ and let $j_1<j_2<\do<j_k$ be a subsequence of $1,2,\do,m$. Let 
$z=i_{j_1}^1i_{j_2}^1\do i_{j_k}^1$ (product in $\fU(\{1\})$); let $\un z$ be the corresponding element of 
$W$. We have $\un z\le w$ for the standard partial order on $W$.}
\nl
We argue by induction on $k$. If $k=0$ the result is obvious. We now assume that $k>0$. We can
assume that $j_1=1$. We set $i=i_1$.
Let $w'=s_{i_2}\do s_{i_m}$. Let $z'=i_{j_2}^1\do i_{j_k}^1$ (product in $\fU(\{1\})$) 
and let $\un z'$ be the corresponding element of $W$. By the induction hypothesis we have $\un z'\le w'$.
Hence $\un z'\le w$.
If $s_i\un z'<\un z'$ then $z=z'$ and $\un z=\un z'$; the inequality $\un z\le w$ follows.
If $s_i\un z'>\un z'$ then $z=s_iz'$ and $\un z=s_i\un z'$. Note that $s_iw<w$. By a standard property
of $\le$ we have $s_iw<w,\un z'\le w\implies s_i\un z'\le w$. Thus $\un z\le w$. This proves (d).

\mpb

Let $\fU(\{1\})^{\le w}$ be the set of elements $z$ of $\fU(\{1\})$ such that the corresponding element
$\un z$ of $W$ satisfies $\un z\le w$. From (c),(d), we see that there is a well defined map 
$$\c_w:\fU_w(\cn)@>>>\fU(\{1\})^{\le w}$$
such that for any $\ii=(i_1,i_2,\do,i_m)\in\co_w$ and any $\aa=(a_1,a_2,\do,a_m)\in\NN^m$ we have
$\c_w(i_1^{a_1}i_2^{a_2}\do i_m^{a_m})=i_{j_1}^1i_{j_2}^1\do i_{j_k}^1$ (product in $U^+(\{1\})$) where 
$j_1<j_2<\do<j_k$ is the sequence consisting of all $j\in[1,m]$ such that $a_j=0$. 

\subhead 10.2\endsubhead
The fibres of $\c_w$ form a partition of $\fU_w(\cn)$ indexed by 
$\fU(\{1\})^{\le w}$. Taking $w=w_I$ we obtain a partition of $\fU_{w_I}(\cn)$ indexed by $W$.
This can be viewed as a partition of $\BB$ (see 0.2) indexed by $W$.

\enddocument